\documentclass[12pt]{article}
\usepackage{amsmath}
\usepackage{amsfonts}
\usepackage{amssymb}
\usepackage{latexsym}
\usepackage{graphicx}
\usepackage{enumerate}
\usepackage[shortlabels]{enumitem}
\usepackage{mathtools}
\usepackage{hyperref}
\usepackage{color}
\usepackage{xcolor}
\usepackage{fullpage}
\usepackage{tikz}

\usepackage{amsmath, amssymb,latexsym}
\usepackage{color}
\usepackage{xcolor}
\usepackage{fullpage}
\usepackage{soul}
\usepackage{float}
\usepackage{comment}

\def\R{\mathbb{R}}
\def\C{\mathbb{C}}

\def\la{\lambda}

\DeclareMathOperator{\diag}{diag}

\DeclareMathOperator{\Rev}{Rev}

\DeclareMathOperator{\rank}{rank}

\newtheorem{theorem}{Theorem}[section]
\newtheorem{proposition}[theorem]{Proposition}

\newtheorem{lemma}[theorem]{Lemma}
\newtheorem{definition}[theorem]{Definition}
\newtheorem{corollary}[theorem]{Corollary}
\newtheorem{remark}[theorem]{{\sc Remark}}
\newtheorem{example}[theorem]{Example}

\title{Generating functions of non-backtracking walks on weighted digraphs: radius of convergence and Ihara's theorem}
\author{Vanni Noferini\thanks{Aalto University, Department of Mathematics and Systems Analysis, P.O. Box 11100, FI-00076, Aalto, Finland. Supported by an Academy of Finland grant (Suomen Akatemian p\"{a}\"{a}t\"{o}s 331240).}
\and 
Mar\'{i}a C. Quintana\thanks{Corresponding author. Aalto University, Department of Mathematics and Systems Analysis, P.O. Box 11100, FI-00076, Aalto, Finland. Supported by an Academy of Finland grant (Suomen Akatemian p\"{a}\"{a}t\"{o}s 331240).}
}

\begin{document}
\maketitle
\begin{abstract} 
It is known that the generating function associated with the enumeration of non-backtracking walks on finite graphs is a rational matrix-valued function of the parameter; such function is also closely related to graph-theoretical results such as Ihara's theorem and the zeta function on graphs. In [P. Grindrod, D. J. Higham, V. Noferini, The deformed graph Laplacian and its application to network centrality analysis, SIAM J. Matrix Anal. Appl. 39(1), 310--341, 2018], the radius of convergence of the generating function was studied for simple (i.e., undirected, unweighted and with no loops) graphs, and shown to depend on the number of cycles in the graph. In this paper, we use technologies from the theory of polynomial and rational matrices to greatly extend these results by studying the radius of convergence of the corresponding generating function for general, possibly directed and/or weighted, graphs. We give an analogous characterization of the radius of convergence for directed (unweighted or weighted) graphs, showing that it depends on the number of cycles in the undirectization of the graph. For weighted graphs, we provide for the first time an exact formula for the radius of convergence, improving a previous result that exhibited a lower bound. Finally, we consider also backtracking-downweighted walks on unweighted digraphs, and we prove a version of Ihara's theorem in that case.
\end{abstract}



\section{Introduction}
Let $G$ be a finite graph, possibly directed and/or weighted, and let $A$ be its adjacency matrix, i.e. $A_{ij}>0$ is the weight of the edge $(i,j)$ in $G$ whereas $A_{ij}=0$ if $(i,j)$ is not an edge in $G$. In this case, it is known that the $(i,j)$ entry of $A^k$ counts the weighted enumeration of distinct walks of length $k$ from node $i$ to node $j$, where the weight of a walk is given by the product of the weigths of its edges. Then the $(i,j)$ entry of the generating function 
\begin{equation}\label{eq:powers}
\sum_{k=0}^\infty t^k A^k
\end{equation} 
 records the weighted sum of all walks from $i$ to $j$, where each walk of length $k$ is counted by its weight times the attenuation factor $t^k$. Katz centrality \cite{Katz} establishes that the importance of a node $i$ can be quantified by the sum over all $j$ of this weighted sum of walks starting from $i$. It is also known then \eqref{eq:powers} converges to $(I-tA)^{-1}$ with radius of convergence $\rho(A)^{-1}$, where $\rho(A)$ denotes the spectral radius of $A$. Then, the centrality of node $i$ is the entry $x_i$ of the vector $x:=(I-tA)^{-1}\pmb{1}$, where $\pmb{1}$ denotes the vector whose entries are all equal to $1$. This leads to the linear system $(I-tA)\pmb{1}=x$.

A walk on a graph is called \emph{non-backtracking} if it does not containg any node subsequence of the form $iji$. The combinatorics of non-backtracking walks has recently received substantial attention in the area of complex network analysis, since it leads to the definition and study of certain \emph{non-backtracking centrality measures} that have some novel features and potential advantages with respect to more traditional ones \cite{AGHN,inprep,nonBT,GHN,pc20,tcter}. Non-backtracking centralities are also linked to the theory of zeta functions on graphs, to Ihara's theorem, and to the so-called Hashimoto (or non-backtracking) matrix: see, for example, \cite{Hashimoto,HST06,ihara,tarfulea} and the references therein. Although not directly relevant to the present article, we mention that the non-backtracking version of the graph Laplacian has also received attention recently \cite{mulas}.

In this paper, we are interested in generalizing \eqref{eq:powers} to non-backtracking walks. We focus on finite graphs, possibly directed, and we also pay attention to the weighted case. Consider first a finite unweighted graph $G$ (directed or undirected) with adjacency matrix $A$, and let $p_k(A)$ be the matrix whose $(i,j)$ element is the number of non-backtracking walks of length $k$ from node $i$ to node $j$. It is known \cite{AGHN} that the generating function 
\begin{equation}\label{eq:genfun}
 \sum_{k=0}^\infty t^k p_k(A)
\end{equation}
converges to a rational matrix-valued function $\Phi(A,t)$ for sufficiently small values of the parameter $t$. For the case of undirected graphs, $\Phi(A,t)$ is linked to the so-called deformed graph Laplacian $M(t):=I-At+(D-I)t^2$, where $D=\diag(\diag(A^2))$: see \cite{Morbidi}. The celebrated Ihara’s Theorem \cite{ihara}, in its original form, relates the determinant of the deformed graph Laplacian with a determinant involving the Hashimoto matrix \cite{Hashimoto} of an undirected graph.

\begin{theorem} [Ihara's Theorem for simple graphs]\label{cor:Ihara} Let $G$ be an unweighted undirected graph of order $n$ with $m$ edges and associated deformed graph Laplacian $M(t)$ and Hashimoto matrix $B$. Then
		\[ \det (I - t B) = (1-t^2)^{m-n} \det M(t).\]
\end{theorem}

As noted in \cite{GHN}, since the generating function in \eqref{eq:genfun} converges to a rational matrix-valued function $\Phi(A,t)$, the radius of convergence $r$ of \eqref{eq:genfun} is given by the smallest pole of $\Phi(A,t)$. For the case of a simple graph, the radius of convergence was characterized in \cite[Theorem 7.3]{GHN} depending on the features of the underlying graph $G$, and in particular on its cycles. In Theorem \ref{thm:old} below, we rephrase that result in an equivalent, but more explicit way that takes into account also \cite[Proposition 4.6 and Theorem 6.3]{GHN} as well as Remark \ref{rem:average} in the present paper.

\begin{theorem}[Theorem 7.3 in \cite{GHN}]\label{thm:old}
    Let $G$ be a simple graph $G$, and let $r$ be the radius of convergence of the generating function in \eqref{eq:genfun}. Then:
   \begin{enumerate}
       \item[\rm (a)] If $G$ has at least one connected component that contains at least two cycles, $r=\mu$, where $\mu \in ]0,1[$ is the smallest eigenvalue of the deformed graph Laplacian $M(t)$;
       \item[\rm (b)] If $G$ contains at least one cycle, but there is no connected component of $G$ that contains more than one cycle, then $r=1$;
       \item[\rm (c)] If $G$ is a forest, $r=\infty$.
   \end{enumerate}
\end{theorem}

A first goal of the present article is to extend Theorem \ref{thm:old} to the case of unweighted directed graphs, that is, to characterize the radius of convergence of \eqref{eq:genfun} with respect to the features of the underlying graph, but without the assumption that the graph is undirected.

It was shown in \cite{AGHN} that the limit of the generating function in \eqref{eq:genfun} is still a rational matrix-valued function also in the case of digraphs (i.e., directed graphs), and that it is linked to a ``directed deformed graph Laplacian" $M(t):=I-At+(D-I)t^2+(A-S)t^3$, where $A$ is the adjacency matrix of the digraph, $D=\diag(\diag(A^2))$ and $S=A \circ A^T$ is the adjacency matrix of the \emph{directed part} of $G$, i.e., the subgraph of $G$ obtained by only keeping the reciprocal edges in $G$. For an undirected graph, $S=A$ and hence one recovers the classical deformed graph Laplacian.

When the edges of the graph (directed or undirected) are assigned positive weights, one can still consider \eqref{eq:genfun} upon replacing the number of NBTW of length $k$ with their weighted sum, where the weight of a walk is product of the weights of its edges. In this setting, formulae for the generating function were obtained in \cite{ryan}. In \cite[Corollary 4.8]{ryan}, a lower bound on the radius of convergence was computed. A second goal of this paper is to prove that this lower bound is in fact equal to the radius of convergence. We also prove a version of Ihara's theorem for weighted digraphs in Theorem \ref{th:det_weights}; see also Remark \ref{rem:notkempton} for why this is result is new and holds in a more general setting than \cite[Theorem 4]{Kempton}.

Finally, rather than fully forbidding backtracking, one could consider a downweighting approach \cite{inprep}. A third goal of our manuscript is to extend the analysis to this case, including a backtrack-downweighted version of Ihara's theorem and a partial extension of Theorem \ref{thm:old}.

The structure of the paper is the following. In Section \ref{sec:preliminary} we develop preliminary graph-theoretical results, including the definitions of undirected part and undirectization of a digraph. In Section \ref{sec:NBTW} we study the radius of convergence of \eqref{eq:genfun} for unweighted digraphs; our main result in this section is Theorem \ref{th:radius}, that in essence shows that Theorem \ref{thm:old} extends to directed graphs provided that the classical deformed graph Laplacian $M(t)$ is replaced by its directed version as in \cite{AGHN}. Furthermore, instead of cycles in $G$, one looks for cycles in the \emph{undirectazion} of $G$, that we define as the graph $H$ obtained from $G$ by making every edge reciprocal (so that $G$ is always a subgraph of $H$). We turn to weighted graphs and weighted digraphs in Section \ref{sec:weighted}. Our main result in this section is Theorem \ref{thm:weightedradius}, that improves \cite[Corollary 4.8]{ryan} by providing an explicit characterization of the radius of convergence of \eqref{eq:genfun} in this setting; moreover, Corollary \ref{thm:weightedradius} is a partial analogue of Theorem \ref{th:radius} in the weighted setting. In Section \ref{sec:BTDW} we analyze the generating function for backtrack-downweighted walks, obtaining a variant of Ihara's theorem in Theorem \ref{th:det} as well as Theorem \ref{th:radius2} which is a partial further extension of Theorem \ref{thm:old} and Theorem \ref{th:radius}. We conclude with some final comments in Section \ref{sec:conc}.

\section{Definitions and preliminary results in graph theory}\label{sec:preliminary}
Let us begin by recalling some basic definitions. All the following concepts are standard, but in our experience the vocabulary of graph theorists appears to have a high variability depending on an individual author's taste, and therefore it is important to declare our own linguistic preferences.

A finite \emph{directed graph}, or digraph, $G=(V(G),E(G))$ consists of a set $V(G)$ of vertices, also called nodes, and a set $E(G) \subseteq V(G) \times V(G)$ of edges. To lighten the notation, we will often drop the dependence on $G$ and write simply $V$ for the set of vertices and $E$ for the set of edges. The order of a graph is its number of vertices $n:=\# V$. It is important to emphasize that the pairs in the set $E$ are ordered, i.e., for all $i,j \in V$, $(i,j)$ and $(j,i)$ represent two distinct edges and $E$ may contain both, only one, or neither. Directed edges are indicated on the graph with an arrow on the edge, and they are sometimes also called arcs. In graph theory, it sometimes possible to have multiple arcs, i.e., an arc $(i,j)$ may appear multiple times on the digraph; however, our definition excludes this case and we do not consider it in this paper. An arc of the form $(i,i)$ is called a \emph{loop} of $G$. If a graph without loops has order $n$ then the number of 
its directed edges $m:=\# E$ satisfies $0 \leq m \leq n^2-n$.

For an unweighted digraph of order $n$, the \emph{adjacency matrix} $A$ of $G$ is an $n\times n$ matrix with $A_{ij}=1$ if $(i,j)\in E$ and $A_{ij}=0$ otherwise. A directed edge $(i,j) \in E$ is called reciprocated if $(j,i) \in E$: in that case the pair of edges is collectively referred to as a reciprocal edge, and we say that $(i,j)$ is the reverse edge of $(j,i)$, and vice versa. An \emph{undirected}  graph is a digraph such that all its edges are reciprocated. Note that a  graph is undirected if and only if its adjacency matrix is symmetric.

A \emph{walk} of length $r$ is a sequence of $r$ directed edges $ (e_1, e_2, …, e_{r})$ for which there is a sequence of vertices $(v_1, v_2, \ldots , v_{r+1})$ such that $e_k=(v_{k},v_{k+1})\in E$ for all $k =1, 2, \ldots , r$. A walk can be either closed (if $v_1=v_{r+1}$) or open (if $v_1 \neq v_{r+1}$). A walk is said to be \emph{backtracking} (BTW) if there is a sequence of the form $(v_k,v_{k+1},v_k) $ in the sequence of vertices, and \emph{nonbacktracking} (NBTW) otherwise. If all the vertices in a walk are distinct, such a walk is called an \emph{open path}; and if all the vertices in a walk are distinct except that the first and the last one coincide, it is called a \emph{closed path}.

By a \emph{path} we mean any walk which is either an open path or a closed path. A digraph is said to be \emph{strongly connected} if there is a path from every vertex to any other vertex. A \emph{strongly connected component} of a digraph is a maximal strongly connected subgraph. If a node does not belong to any strongly connected component, it is called a \emph{single node}.

We shall also consider weighted (di)graphs. These are pairs $(G,\omega)$ where $G=(V,E)$ is a digraph and $\omega: E \rightarrow ]0,+\infty[$ is a function that assigns a positive weight $\omega(e)>0$ to each edge $e \in E$. The weight of a walk $(e_1,\dots,e_r)$ is defined multiplicatively as $\prod_{i=1}^r \omega(e_i)$; note that an edge is allowed to possibly appear more than once in the sequence $e_1,\dots,e_r$, and hence the weight of a walk can equivalently be written as $\prod_{j=1}^\ell (\omega(e_{i_j}))^{m_j}$ where now $e_{i_1},\dots,e_{i_\ell}$ are the distinct edges in the walk and $m_j$ is the number of times that edge $e_{i_j}$ is travelled during the walk. The definition of adjacency matrix extends to weighted digraphs by setting $A_{ij}=\omega(e)$ if $e=(i,j) \in E$ and $A_{ij}=0$ otherwise. We observe that it is sometimes convenient to see unweighted digraphs as the special case of weighted digraphs for which the weight function is $\omega(e)=1$ for all $e \in E$.

\subsection{Undirected part and undirectization of a digraph}

In this subsection, we define two undirected graphs that we associate with each digraph $G$. The first one is the \emph{undirected part} of $G$, denoted by $G_U(G)$ and obtained by keeping only reciprocated edges in $G$, i.e., by removing edges whose reverse is not present. The second one is the \emph{undirectization} of $G$, denoted by $H(G)$ and obtained by making every edge in $G$ reciprocated, i.e., by adding reverse edges where they are missing. When there are no risks of ambiguity, for the sake of simplicity we may write $G_U$ and $H$ in lieu of $G_U(G)$ and $H(G)$.

\begin{definition}\label{def:undirected}
	Given a digraph $G=(V,E)$, the \emph{undirected part} of $G$ is the graph $G_U(G)=(V,E_U)$ where $E_U = \{ (i,j) \in E : (j,i) \in E \}$. Moreover, the \emph{undirectization} of $G$ is the graph $H(G)=(V,E_H)$ where $E_H=\{ (i,j) : (i,j) \in E \ \mathrm{or} \ (j,i) \in E\}$.
\end{definition}

\begin{example}\label{example1}

Let $G=(V,E)$ with $V=[4]=\{1,2,3,4\}$ and $E=\{ (1,2),(2,1),(2,3),(3,4),\allowbreak (4,2)\}$. Then, the undirected part of $G$ is $G_U=(V,E_U)$ with $E_U=\{(1,2),(2,1) \}$ whereas the undirectization of $G$ is $H=(V,E_H)$ with $E_H=\{ (1,2),(2,1),(2,3),(2,4),(3,2),(3,4),(4,2),\allowbreak (4,3) \}$. 

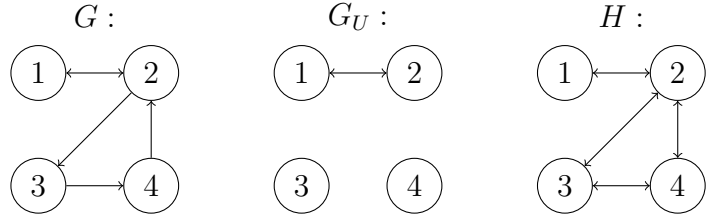
\begin{figure}[h]
		\begin{center}
			\begin{tikzpicture}[scale=0.5]
			
			\node[draw,circle] (1) at (-7,3) {$1$};
			\node[draw,circle] (2) at (-4,3) {$2$};
			\node[draw,circle] (3) at (-7,0) {$3$};
			\node[draw,circle] (4) at (-4,0) {$4$};
			\draw[<->] (1) -- (2) ;
			\draw[->] (2) -- (3);
			\draw[->] (3) -- (4);
			\draw[->] (4) -- (2);
			\draw (-5.5,4.5)  node {$G:$} ;

			\node[draw,circle] (11) at (0,3) {$1$};
			\node[draw,circle] (22) at (3,3) {$2$};
			\node[draw,circle] (33) at (0,0) {$3$};
			\node[draw,circle] (44) at (3,0) {$4$};
			\draw[<->] (11) -- (22);
			\draw (1.5,4.5)  node {$G_U:$} ;
			
			\node[draw,circle] (111) at (7,3) {$1$};
			\node[draw,circle] (222) at (10,3) {$2$};
			\node[draw,circle] (333) at (7,0) {$3$};
			\node[draw,circle] (444) at (10,0) {$4$};
			\draw[<->] (111) -- (222);
			\draw[<->] (222) -- (333);
			\draw[<->] (333) -- (444);
			\draw[<->] (444) -- (222);
			\draw (8.5,4.5)  node {$H:$} ;

			\end{tikzpicture}\end{center}
		
		\caption{An unweighted digraph $G$, its undirected part $G_U$ and its undirectization $H$.}
		\label{fig:undirect}
		
		\end{figure}

The (di)graphs $G,G_U$ and $H$ have adjacency matrices, respectively, $A,S=A \circ A^T$ and $A+A^T-S$:

\[ A=\begin{bmatrix}
    0&1&0&0\\
    1&0&1&0\\
    0&0&0&1\\
    0&1&0&0
\end{bmatrix}, S=\begin{bmatrix}
    0&1&0&0\\
    1&0&0&0\\
    0&0&0&0\\
    0&0&0&0
\end{bmatrix}, A+A^T-S=\begin{bmatrix}
    0&1&0&0\\
    1&0&1&1\\
    0&1&0&1\\
    0&1&1&0
\end{bmatrix}.      \]

\end{example}

It is clear by Definition \ref{def:undirected} that $G_U$ is a subgraph of $G$ which is in turn a subgraph of $H$. Moreover, as it is illustrated by Example \ref{example1}, if $G$ has adjacency matrix $A$ then $G_U$ has adjacency matrix $S:=A\circ A^T$, while $H$ as adjacency matrix $A+A^T-S$. Let $D:=\diag (\diag A^2)$ be the diagonal matrix whose diagonal entry $d_{ii}$ is the number of reciprocal edges departing from node $i$. The matrix $D$ is called the \emph{degree matrix} of $G$. Note that $D$ and $S$ satisfy the relation $D=\diag (\diag S^2)$. In other words, $D$ is also the degree matrix of $G_U$. Instead, the degree matrix of $H$ is $D_H := \diag (\diag\,(A+A^T-S)^2) \geq D$. Moreover, $D_H=D$ if and only if $G_U=G=H$.

\subsection{Directed cycles}

A \emph{directed cycle} in a digraph is a directed version of a cycle in a undirected graph, with all the edges oriented in a coherent direction. 
Namely, a directed cycle of length $\ell \geq 3$ is a closed path of length $\ell \geq 3$, modulo any cyclic permutation on its vertices. (In other words, for example, the three vertex sequences $(1,2,3,1)$, $(2,3,1,2)$ and $(3,1,2,3)$ all represent the same cycle). A \emph{directed acyclic graph} is a digraph with no directed cycles. A directed acyclic graph $G$ is said to be an \emph{undirected tree} if $G=G_U,$ where $G_U$ is the undirected part of $G$, and $G_U$ is a tree.  {Note that, if $G$ is strongly connected, then its undirectization $H$ is connected. In Lemma \ref{lem:undirectedtree}, we prove that a strongly connected digraph $G$ is an undirected tree if and only if its undirectization $H$ is a tree.

\begin{lemma}\label{lem:undirectedtree}
    Let $G$ be a strongly connected digraph, and let $H$ be the undirectization of $G$. Then, $H$ is a tree if and only if $G$ is an undirected tree, and in this case $H=G$.
\end{lemma}

\begin{proof} If $G$ is an undirected tree, it is obvious by definition that $H=G$ is a tree. For the other implication, assume that $H$ is a tree. If $G$ is not an undirected tree then there are vertices $i,j$ such that $(i,j) \in E(G)$ but $(j,i) \not \in E(G)$; however, $G$ is strongly connected, so there is a path $(j,i_1,\dots,i_{\ell-1},i)$ on $G$ of length $\ell > 1$ (as there must be a path but there is no path of length $1$ by assumption). But this contradicts that $H$ is a tree because it implies that $(j,i_1,\dots,i_{\ell-1},i,j)$ is a cycle of length $\ell+1 \geq 3$ on $H$. Finally, having established that $G$ is an undirected tree, it is clear that $H=G$.
\end{proof}
}

If $\Gamma=(i_0,i_1,\dots,i_\ell=i_0)$ is a directed cycle of length $\ell$, the set of vertices of $\Gamma$ is the (unordered) set $\{i_0,i_1,\dots,i_{\ell-1}\}.$ Note that the set of vertices is uniquely determined by a directed cycle, even though the sequence of vertices that defines a directed cycle $\Gamma$ is not unique (as one might apply a cyclic permutation to the vertices). {However, the set of vertices does not uniquely identify a directed cycle.} We say that two directed cycles $\Gamma_1$ and $\Gamma_2$ are the reverse of each other if (up to some cyclic permutations) $\Gamma_1=(i_0,i_1,\dots,i_{\ell-1},i_\ell=i_0)$ and $\Gamma_2=(i_\ell=i_0,i_{\ell-1},\dots,i_1,i_0)$. We say that two directed cycles in a graph have at least one different vertex if the set of vertices, say, $V_1$ and $V_2$ in the two cycles are distinct. Note that this is compatible with either $V_1\subset V_2$ or $V_2 \subset V_1$. In particular, cycles of distinct lengths always have at least one different vertex. Directed graphs may have distinct directed cycles that do not have different vertices, such as two cycles whose edges are all reciprocal, and that hence can be travelled in both orientations. It may also happen that the vertices are travelled in a permuted order, which is neither a cyclic permutation nor an inversion. For example, $(1,2,3,4,1)$ and $(1,3,4,2,1)$ are distinct cycles with the same set of vertices. If two directed cycles in a digraph $G$ have the same vertices and are not each other's reverse then we can guarantee that $G$ has at least two directed cycles with at least one different vertex as the next lemma shows. For instance, in the example above, one can easily show that necessarily the graph contains also the cycles $(1,3,4,1)$ and $(2,3,4,2)$. The general case is analyzed in Lemma \ref{lem:eitherdifferentorreverse} below.

\begin{lemma}\label{lem:eitherdifferentorreverse}
	Let $G$ be a digraph with at least two distinct directed cycles. Then either the two cycles are one the reverse of the other or $G$ has at least two cycles with at least one different vertex.
\end{lemma}

\begin{proof}
	Let $C_1,C_2$ be two distinct directed cycles in $G$. If the sets of their vertices are distinct then there is nothing to prove. Thus, assume that the vertices of $C_1$ and $C_2$ are the same and that $C_1$ and $C_2$ are not each other's reverse. Note that, in such a case, the length of these cycles must necessarily be $\ell \geq 4$. Then, there is a starting vertex $v_1$ and permutation $\sigma$ on $\{2,\dots,\ell\}$, not equal to the identity, such that $\sigma(\ell)\neq 2$, and such that the two cycles are
	\[C_1 := v_1 \rightarrow v_2 \rightarrow \dots \rightarrow v_\ell \rightarrow v_1 \]
	and
	\[C_2:= v_{1} \rightarrow v_{\sigma(2)} \rightarrow \dots \rightarrow v_{\sigma(\ell)} \rightarrow v_{1},\]
	respectively. Let $j = \min \{ i \, : \, i \neq \sigma(i) \}$. Then, $\sigma(k)=j$ for some $k>j$ and therefore
	\[ C_3:=v_1 \rightarrow v_2 \rightarrow \dots \rightarrow v_j \rightarrow v_{\sigma(k+1)} \rightarrow \dots \rightarrow v_{\sigma(\ell)} \rightarrow v_1\]
	is a closed path of length $t:=\ell-k+j$. (If $k=\ell$, we agree by convention that $\sigma(k+1)=1$.) Since $k>j$ then $t<\ell$. On the other hand, $t>2$ since $\sigma(\ell)\neq 2$. We conclude that $C_3$ is a cycle of length strictly less than $\ell$, and hence $C_1$ and $C_3$ are two cycles with at least one different vertex.
\end{proof}

For an undirected graph $G$, if $\Gamma$ is a directed cycle of $G$ then the reverse of $\Gamma$ is also a directed cycle of $G$. A \emph{cycle} on an undirected graph is a closed path of length $\ell \geq 3$ up to cyclic permutations \emph{or inversions}. For example, on an undirected graph, $(1,2,3,1)$ and $(1,3,2,1)$ represent the same cycle (although they represent distinct directed cycles).

\begin{lemma}\label{lem:HimpliesG} Let $G$ be a strongly connected digraph and let $H$ be the undirectization of $G$. If $H$ has at least $2$ cycles with at least one different vertex, then $G$ has at least $2$ directed cycles with at least one different vertex.
	
\end{lemma}

\begin{proof}
	Each of the two cycles in $H$ induce a directed cycle in $G$ if and only if there is a corresponding set of edges in $G$ all with coherent orientation. If this happens for both cycles, then we have identified two cycles in $G$ with at least one different vertex. Thus we assume, for instance, that two edges in the set of edges corresponding to the first cycle are not reciprocated and do not have coherent orientation. Then we have in $G$ two edges of the form
	$$v_{i-1}\leftarrow v_{i} \rightarrow v_{i+1},$$
	and moreover they are both not reciprocated.
	Since $G$ is strongly connected, we can construct a path from $v_{i-1}$ to $v_{i}$ (that, together with the edge $v_{i}\rightarrow v_{i-1}$, forms a directed cycle on $G$, that we call $C_-$) and we can construct a path from $v_{i+1}$ to $v_{i}$ (that, together with the edge $v_{i} \rightarrow v_{i+1}$, forms a directed cycle on $G$, that we call $C_+$). It is clear that $C_+$ and $C_-$ are distinct, as $C_- \not \ni v_i \rightarrow v_{i+1} \in C_+$. In addition, $C_+$ and $C_-$ are not the reverse of each other since $v_i \rightarrow v_{i+1} \in C_+$ but $v_{i+1} \rightarrow v_i \not \in E(G) \Rightarrow v_{i+1} \rightarrow v_i \not \in C_-$. Then, we have two distinct directed cycles in $G$ that are not one the reverse of the other and, thus, Lemma \ref{lem:eitherdifferentorreverse} concludes the proof.
\end{proof}

\section{Unweighted digraphs: NBTW walks and the directed deformed graph Laplacian}\label{sec:NBTW}

 Let $G$ be an unweighted digraph without loops. We consider $p_k(A)$ the matrix whose entry $(p_k(A))_{ij}$ counts the number of NBTWs of length $k$ from node $i$ to node $j$. It was proved in \cite{AGHN} that the associated generating function has the following expression:
$$\displaystyle\sum_{k=0}^{\infty} t^k p_{k}(A)=(1-t^2)(I - A t + (D-I) t^2 + (A-S) t^3)^{-1},$$
where $S=A\circ A^T$ and $D=\diag (\diag A^2)$ is the diagonal matrix counting reciprocal edges from/to each node in $G$. If the graph is undirected then the expression above for the generating function is simplified as $S=A$. The polynomial matrix \begin{equation}\label{eq:dgl}
    M(t):= I - A t + (D-I) t^2 + (A-S) t^3
\end{equation} is called the \emph{directed deformed graph Laplacian}\footnote{In the literature, the name deformed graph Laplacian is only used in the undirected case $A=S$. In this paper, however, we find it convenient to also give a name to the expression \eqref{eq:dgl}} of $G$. Note that, in the special case where $G$ is undirected, then $S=A$ and hence we recover the same expression studied in \cite{GHN} and mentioned in the introduction. We will characterize the radius of convergence of the generating function and the spectrum of $M(t)$ with respect to the underlying digraph. We now briefly recall some relevant definitions about the spectral theory of polynomial matrices, see e.g. \cite{indexsum,GLR82, GHN, GT, nnt} and the references therein.

\begin{remark}
    We observe how the directed deformed graph Laplacian of $G$ as in \eqref{eq:dgl} compares with the deformed graph Laplacian of $G_U$, the undirected part of $G$, and of $H$, the undirectization of $G$.

The former has deformed graph Laplacian
\begin{equation}\label{eq:1}
	M_U(t) = I - S t + (D-I) t^2 =  M(t) - (t^3-t) (A-S).
	\end{equation}
    while the latter has deformed graph Laplacian
\[ M_H(t)=I-(A+A^T-S) t + (D_H - I)t^2 = M(t) - (A^T-S) t  + (D_H-D) t^2. \]

\end{remark}

\subsection{Background material on polynomial matrices}

A polynomial matrix $M(t) \in \C[t]^{n \times n}$ is called regular if $\det M(t) \not \equiv 0$. A finite number $\lambda \in \C$ is an \emph{eigenvalue} of the regular $M(t)$ if it is a root of $\det M(t)$. The \emph{algebraic multiplicity} of an eigenvalue $\lambda$ is its multiplicity as a root of the scalar polynomial $\det M(t)$. Infinity is an eigenvalue of a polynomial matrix $M(t)$ if $0$ is an eigenvalue for the reversal polynomial matrix $\Rev M(t):=t^k M(\frac{1}{t})$, and the algebraic multiplicity of $\infty$ is defined as that of $0$ for $\Rev M(t)$.

For every polynomial matrix $M(t) \in \C[t]^{n \times n}$ having rank $r$ over the field $\C(t)$, there exist two unimodular (that is, invertible over the ring of polynomials) polynomial matrices $U(t)$ and $V(t)$ such that $U(t)M(t)V(t)=S(t)$ where $S(t)$ is diagonal and its diagonal elements $S_{ii}(t)$, for $i=1,\ldots,r,$ are monic polynomials that form a divisibility chain, that is, $S_{jj}(t)$ divides $S_{j+1,j+1}(t)$, for $j=1,\dots,r-1$; and $S_{ii}(t)=0$ for $i=r+1,\dots,n$. The matrix $S(t)$ is called the \emph{Smith form} of $M(t)$, and is uniquely determined by $M(t)$. The diagonal elements of the Smith form of $M(t)$ are called the \emph{invariant polynomials} of $M(t)$. It is clear that, when $M(t)$ is regular the diagonal elements of $S(\la)$ are all nonzero and their product is (up to multiplication by a {nonzero} constant) equal to $\det M(t)$. The \emph{geometric multiplicity} of a finite eigenvalue $\lambda$ is the number of invariant polynomials of which $\lambda$ is a root; and the \emph{partial multiplicities} of $\lambda$ are the (positive) multiplicities of $\lambda$ as a root of each invariant polynomial. The partial multiplicities and the geometric multiplicity of $\infty$ for the polynomial matrix $M(t)$ of degree $k$ are defined as, respectively, the partial multiplicities and the geometric multiplicity of $0$ for the reversal polynomial matrix $\Rev M(t)$. 

\begin{example} 
Consider the following polynomial matrix:
\[ M(t) = \begin{bmatrix}
    1 & 0 & 0\\
    0 & t & 0\\
    0 & 0 & t^3-2t^2+t
\end{bmatrix},\]
which is already in Smith form. It is easy to see that the partial multiplicities of the eigenvalue $0$ are $1$ and $1$, and hence the algebraic multiplicity of $0$ is $2$ and the geometric multiplicity is $2$. Similarly, the partial multiplicity of $1$ is $2$. Computing the partial multiplicities of $\infty$ requires a little more work, but the readers can convince themselves that they are $2$ and $3$. The relation $(1+1+2+2+3)=3 \cdot 3$ is not coincidental: for a regular polynomial matrix of degree $k$ and size $n$, the sum of the partial multiplicities of all the (finite and infinite) eigenvalues is always equal to $kn$. This fact is known as the Index Sum Theorem \cite{indexsum}.
    
\end{example}

\subsection{Eigenvalues of the directed deformed graph Laplacian: Algebraic, geometric and partial multiplicities}
In this subsection, we discuss some properties of the eigenvalues of the directed deformed graph Laplacian $M(t)$ in \eqref{eq:dgl} and their multiplicities. We first note that, given an unweighted digraph with no loops $G=(V,E)$, the size of $M(t)$ is $n=\#V$ whereas its degree $k$ is $k=3$ if and only if $G$ is not undirected, and $k=2$ if and only if $G$ is undirected. As the undirected case was already studied in \cite{GHN}, we will henceforth assume that $k=3$.

Recall that a nonnegative square matrix is said to be reducible if there exists a permutation matrix $P$ such that \[ P^TA P= \begin{bmatrix}
A_1 & E\\
0 & A_2
\end{bmatrix} ,\] where $A_1$ and $A_2$ are square and have positive size. Otherwise $A$ is said to be irreducible. It is known that an adjacency matrix is reducible if and only if its associated digraph is not strongly connected. In addition, if the adjacency matrix $A$ of a digraph $G$ is reducible then it can be written in the form
  \[ A = \begin{bmatrix}
  A_1 & E_{12} & \cdots & E_{1k}\\
  0 & A_2 & \cdots & E_{2k} \\
  \vdots & \vdots & \ddots & \vdots \\
  0& 0 & \cdots & A_k
  \end{bmatrix} ,\] 
  up to row/column permutations, where $A_1, A_2, \cdots , A_k$ are { either square irreducible matrices of positive sizes or $1 \times 1$ zero matrices}. Then the digraphs $G_i$ associated with $A_i$ are either strongly connected subgraphs {or single nodes} of $G$. If a digraph $G$ is not strongly connected, without loss of generality, {for example by separating one strongly connected component/single node from the rest of the graph,} we may therefore assume that its adjacency matrix $A$ can be written in the form  
\[ A= \begin{bmatrix}
A_1 & E\\
0 & A_2
\end{bmatrix} .\] Then

\[A^2 = \begin{bmatrix}
A_1^2 & \star \\
0 & A_2^2
\end{bmatrix}\quad  \text{ and }\quad  S = \begin{bmatrix}
S_1 & 0\\
0 & S_2
\end{bmatrix},\]
where $\star$ denotes a block whose exact nature is irrelevant and $S_i = A_i \circ A_i^T$. Denoting by $M_i(t)$ the directed deformed graph Laplacian of the digraph $G_i$ associated with $A_i$, we deduce that the directed deformed graph Laplacian of $G$ can be written as
\begin{equation}\label{eq:components}
    M(t) = \begin{bmatrix}
M_1(t) & (t^3-t) E\\
0 & M_2(t) 
\end{bmatrix} .
\end{equation}

{The discussion above implies the following Proposition \ref{prop:algebraic} on the algebraic multiplicities of the eigenvalues of $M(t)$.

\begin{proposition}\label{prop:algebraic}
    Let $M(t)$ be the directed deformed graph Laplacian of a digraph $G$, and let $M_{i}(t)$ be the directed deformed graph Laplacian of each strongly connected or single node component $G_i$ of $G$. Then the algebraic multiplicities of a finite eigenvalue $\lambda$ of $M(t)$ is the sum of the algebraic multiplicities of $\lambda$ for each $M_{i}(t)$.
\end{proposition}
\begin{proof} The results follows by noting that, if $G$ is not strongly connected, $M(t)$ can be recursively written as in \eqref{eq:components}.   
\end{proof}

\begin{remark}\rm \label{rem:algebraic}
 It follows by Proposition \ref{prop:algebraic} that, for the purpose of discussing the \emph{algebraic} multiplicities of eigenvalues of the directed deformed graph Laplacian $M(t)$, there is no loss of generality in assuming that the graph $G$ is strongly connected, for if not we can study the algebraic multiplicities of the eigenvalues of the directed deformed graph Laplacians of each strongly connected or single node component $G_i$. 
\end{remark}}

\begin{remark}\rm \label{rem:geometric1}
For the eigenvalues $\pm 1$, the geometric multiplicities is also the sum of the geometric multiplicities for each {directed} deformed graph Laplacians $M_{i}(t)$ of each strongly connected {or single node} component since $M(\pm 1)= M_1(\pm 1) \oplus M_2(\pm 1)$. Indeed, note that, almost immediately by its definition, the geometric multiplicity of a finite eigenvalue $\lambda$ is equal to the nullity of $M(\lambda)$.
\end{remark}
 In addition, we have the following result.
\begin{proposition}
Let $M(t)$ be the {directed} deformed graph Laplacian of a digraph $G$. Then, $1$ is always an eigenvalue, and its geometric multiplicity is the number of connected components \footnote{For an undirected graph, an isolated node is conventionally considered a connected (bipartite) component of $G_U$, the undirected part of $G$.} Moreover, $-1$ is an eigenvalue if and only if $G_U$ has at least one bipartite connected component. The geometric multiplicity of $-1$ is equal to the number of such components.
\end{proposition}

\begin{proof}Denoting by $M(t)$ and $M_U(t)$ the {directed} deformed graph Laplacians of $G$ and $G_U$, respectively, we have by \eqref{eq:1} that $M(t) = M_U(t) + (t^3-t) (A-S).$	 Let us now denote by $L:=D-S$ and $Q:=D+S$, respectively, the graph Laplacian and the signless Laplacian of $G_U$, the undirected part of $G$. Then, by \eqref{eq:1}, $M(1) = M_U(1) = L$. Let $e$ be the vector of all ones. Then $Le=0$ and the nullity of $L$, i.e., the geometric multiplicity of $1$, is the number of connected component of $G_U$. Moreover, we have that $M(-1)=M_U(-1)=Q$ by \eqref{eq:1}, and $Q$ is singular if and only if $G_U$ has at least one bipartite component. In addition, the geometric multiplicity is the number of such bipartite components. 	
\end{proof}

{ It will also be useful in following sections to characterize when $1$ is a defective eigenvalue of the directed deformed graph Laplacian}. For that, we set 
\begin{equation}\label{def:total}
d := e^T A e \text{ , i.e., the total number of directed edges in } G,
\end{equation}
also equal to both the sum of out-degrees and the sum of in-degrees of $G$, and 
\begin{equation}\label{def:totalrecip}
d_U := e^T S e \text{ , i.e., the total number of reciprocated edges in } G,
\end{equation}
also equal to {twice} the sum of all the degrees in the undirected part $G_U$ of $G$.

\begin{remark}\label{rem:average} \rm
	Note that, if $G$ is strongly connected, then its undirectization $H$ is connected and $(2d-d_U)/n$ is equal to the average degree of $H$, where $n$ is the order of $G$. In addition, we have that
	\begin{itemize}
		\item $2d-d_U > 2 n \Leftrightarrow$ the average degree of $H$ is higher than $2$ $\Leftrightarrow H$ has at least two cycles;
		\item $2d-d_U = 2 n \Leftrightarrow$ the average degree of $H$ is equal to $2 \Leftrightarrow H$ has precisely one cycle;
		\item $2d-d_U < 2 n \Leftrightarrow$ the average degree of $H$ is smaller than $2 \Leftrightarrow H$ is a tree.
	\end{itemize}
\end{remark}

\begin{proposition}\label{prop:defective} Let $G$ be a digraph of order $n$ with $d,d_U$ defined as above, and let $M(t)$ be the directed deformed graph Laplacian associated with $G$. Then $1$ is a defective eigenvalue of $M(t)$ if and only if $2 d - d_U = 2n$.
\end{proposition}
\begin{proof}
	$1$ is a defective eigenvalue of $M(\lambda)$ if and only if there exists a vector $v$ such that $M(1)v+M'(1)e=0$. Since $M'(1)=2L+ 2A - 2 I - S$ then $M(1)v+M'(1)e=0$ for some $v$ if and only if  $0 = e^T (2A - 2 I - S) e  = 2d -2n - d_U$.
\end{proof}

By Remark \ref{rem:average} and Proposition \ref{prop:defective}, we obtain the following corollary.

\begin{corollary}\label{cor:defective} Let $G$ be a strongly connected digraph and let $M(t)$ be its directed deformed graph Laplacian. Let $H$ be the undirectization of $G$. Then $1$ is a defective eigenvalue of $M(t)$ if and only if $H$ has precisely one cycle.
\end{corollary}

{Finally, we study the geometric multiplicities of those eigenvalues $\lambda \neq \pm 1$.}

\begin{proposition} Let $G$ be a digraph and let $M(t)$ be its directed deformed graph Laplacian.
 If $\lambda \neq \pm 1$ is a finite eigenvalue of only one directed deformed graph Laplacian $M_{i}(t)$ associated with a strongly connected component $G_i$ of $G$, then its geometric multiplicity as an eigenvalue of $M(t)$ and as an eigenvalue of $M_{i}(t)$ is the same.
\end{proposition}
\begin{proof}
	 Suppose first that $G$ has two strongly connected components (the case where $G$ is strongly connected is trivial). If $\lambda$ is an eigenvalue only for $M_1(t)$ then $M_1(\lambda) U=0$ if and only if
	\[  \begin{bmatrix}
	M_1(\lambda) & (\lambda^3-\lambda) E\\
	0 & M_2(\lambda)
	\end{bmatrix} \begin{bmatrix}
	U\\
	0
	\end{bmatrix} = 0\]
	for any full column rank matrix $U$, which columns would be eigenvectors associated with $\lambda$. On the other hand, if $\lambda$ is an eigenvalue only for $M_2(t)$ then $M_2(\lambda)V=0$ if and only if
	\[  \begin{bmatrix}
	M_1(\lambda) & (\lambda^3-\lambda) E\\
	0 & M_2(\lambda)
	\end{bmatrix} \begin{bmatrix}
	(\lambda-\lambda^3) M_1(\lambda)^{-1} E V\\
	V
	\end{bmatrix} = 0,\]
	for any full column rank matrix $V$.

 The case of more than two strongly connected components can be dealt with by iteratively applying the argument above.
\end{proof}

However, if $\lambda \neq \pm 1$ is a finite eigenvalue of both {directed} deformed graph Laplacians $M_1(t)$ and $M_2(t)$ of each strongly connected component, then its geometric multiplicity as an eigenvalue of $M(t)$ may not be equal to the sum of the geometric multiplicities (see Example \ref{ex_multiplicities}). { Note also that Example \ref{ex_multiplicities} illustrates that} the partial multiplicities may be different than the union of the partial multiplicities of each strongly connected component even for $\pm 1$. 

\begin{example}\label{ex_multiplicities}
Consider the digraph with adjaceny matrix

\[ A = \begin{bmatrix}
B & I\\
0 & B
\end{bmatrix}, \qquad B = \begin{bmatrix}
0 & 1 & 1\\
1 & 0 & 1\\
1 & 1& 0
\end{bmatrix}.\]

Then the two strongly connected components, associated with $B$, have deformed graph Laplacian $M_B(t)=I - t B + I t^2 + 0 t^3$ whose Smith form is {$\diag(1 , x^2+x+1 , (x^2+x+1)(x-1)^2)$}. The whole graph has deformed graph Laplacian

\[ M_A(t) = \begin{bmatrix}
M_B(t) & (t^3-t) I \\
0 & M_B(t)
\end{bmatrix}   \]
whose Smith form is
 { $\diag( I_4 , (x-1)(x^2+x+1)^2 , (x-1)^3 (x^2+x+1)^2).$}
So the partial multiplicities of the eigenvalue $1$ are $0,0,2$ for each connected component, but $0,0,0,0,1,3$ for the whole graph. Note that the same example also shows that geometric multiplicities of eigenvalues not equal to $\pm 1$ may differ from the sum of the geometric multiplicities of the components, and that the partial multiplicities of eigenvalues not equal to $\pm$ 1 may also differ from the union of the partial multiplicities of the components.

\end{example}

\subsubsection{Pruning reciprocal leaves}

We say that a vertex $i \in V$ is a \emph{reciprocal leaf} of $G$ (and $G_U$) if it only adjacent to another vertex $j \in V$ and both $(i,j) \in E$ and $(j,i) \in E$. We show in Theorem \ref{thm:leaves} that pruning reciprocal leaves does not change geometric and algebraic multiplicity of finite eigenvalues. {Theorem \ref{thm:leaves} generalizes \cite[Theorem 5.4]{AGHN}, which however only analyzed the algebraic (but not the geometric) multiplicities.}

\begin{theorem}\label{thm:leaves}
Let $G$ be a digraph with a reciprocal leaf $i$, and let $\widetilde{G}$ be the graph obtained from $G$ by removing $i$ and the unique reciprocal edge connecting $i$ with the rest of $G$. Let $M(t)$ and $\widetilde{M}(t)$ be the directed deformed graph Laplacians associated with $G$ and $\widetilde{G}$, respectively. Then, the finite eigenvalues of $M(t)$ and $\widetilde{M}(t)$ are the same and have the same algebraic and geometric multiplicities.
\end{theorem}

\begin{proof}
There is no loss of generality in assuming that $i=1$ and that the node is adjacent to node $2$. Then,
\[ 1 \oplus \widetilde{M}(t) = M(t) -\left( \begin{bmatrix}
0 & -t\\
-t & t^2
\end{bmatrix} \oplus 0 \right) \]
implying by the matrix determinant lemma
\[ \frac{\det \widetilde{M}(t)}{\det M(t) } = \det \left( I_2 - N(t) \begin{bmatrix}
0 & -t\\
-t & t^2
\end{bmatrix} \right),    \]
where $N(t)$ is the top $2 \times 2$ block of $M(t)^{-1}$. Each entry $[N(t)]_{ij}$ is equal to $(1-t^2)^{-1}$ times the weighted sum of nonbacktracking walks (on $G$) from node to $i$ to node $j$, for $i,j\in\{1,2\}$. Let $1+f(t)$ be such a sum for $i=j=2$. Since $1$ is a reciprocal leaf adjacent to $2$, it is easy to see that
\[ N(t) = \frac{1}{1-t^2} \begin{bmatrix}
1+t^2 f(t) &t + t f(t) \\
t + t f(t) & 1+f(t)
\end{bmatrix}.\]
After some straigthforward manipulations we obtain that $ \dfrac{\det \widetilde{M}(t)}{\det M(t) } = 1, $ which proves the statement for algebraic multiplicities.

Concerning the geometric multiplicity, suppose that $M(\lambda) V = 0$ for some $\lambda \in\C$ and some full column rank $V$, and partition \[ V= \begin{bmatrix}
v_1 \\v_2 \\ \widetilde{V}
\end{bmatrix}= \begin{bmatrix}
v_1 \\W
\end{bmatrix} \] where $v_1,v_2$ are row vectors. Then,
\[ M(\lambda) V = 0 \Leftrightarrow \begin{bmatrix}
v_1\\
\widetilde{M}(\lambda)W 
\end{bmatrix} = \begin{bmatrix}
\lambda v_2 \\
\lambda v_1 - \lambda^2 v_2 \\
0
\end{bmatrix} \Leftrightarrow v_1 = \lambda v_2  \ \mathrm{and} \ \widetilde{M}(\lambda)W=0. \]
Moreover, 
\[ W = V \begin{bmatrix}
0\\
I
\end{bmatrix}  \Rightarrow \rank W = \rank V.\]
Taking $V$ to be the matrix of all the eigenvectors of $M(t)$ associated with $\lambda$ and $W$ the matrix of all the eigenvectors of $\widetilde{M}(t)$ associated with $\lambda$ yields the statement.
\end{proof}

%

\subsection{Non-$k$-cycling matrices}\label{sub:nonkcycling}

In this subsection we consider the non-$k$-cycling matrices $P_k$ in \cite{nonBT} associated with a graph $G$, and we give a generalization of a theorem for undirected graphs in \cite{nonBT} to directed graphs. Let $G$ be an unweighted digraph without loops. 
We denote by $B\in\mathbb{R}^{m\times m}$ the {\em Hashimoto matrix} \cite{Hashimoto} of $G$ or the {\em non-backtracking edge-matrix}. The Hashimoto matrix is the adjacency matrix of a graph of $m$ nodes, each corresponding to an edge in $G$, obtained by connecting two of the nodes if and only if the corresponding two edges form a NBTW of length $2$ in $G$. This definition is generalised as follows.

\begin{definition}[Non-$k$-cycling matrix] Let $G$ be a digraph and let $k\geq 2$. The non-$k$-cycling matrix $P_k\in\mathbb{R}^{p\times p}$ of $G$ is a matrix whose size $p$ is the number of open paths of length $k-1$ in $G$. Then, the each entry $(P_k)_{ij}$ is associated with two open paths $i=(i_1,\ldots,i_{k})$ and $j=(j_1,\ldots,j_{k})$ of length $k-1\geq 1$ and takes the value $1$ if $G$ has an open path of length $k$ from node $i_1$ to node $j_k$ which superposes with $i$ on its last $k-1$ steps and with $j$ on its first $k-1$ steps. That is,
	
    $$(P_k)_{ij}= \left\{ \begin{array}{lcc}
	1 &   \text{ if } i_1\neq j_k \text{ and } i_r=j_{r-1}\text{ for } r=2,\ldots,k     \\
	0 &  otherwise 
	\end{array}
	\right.$$

\end{definition} 

For $k=2$, the non-$k$-cycling matrix is the Hashimoto matrix, i.e., $P_2=B$. For $k=1$, the non-$k$-cycling matrix is defined as the adjacency matrix of $G$, i.e., $P_1=A$. 

\begin{theorem} \label{thm:directedPk} Let $G$ be a strongly connected digraph and let $P_k$ be the non-$k$-cycling matrix of $G$ for $k\geq 2$. If $G$ has at least $2$ cycles of length greater than $k$ with at least one different vertex then $\rho (P_k) > 1$.
\end{theorem}

\begin{proof} We will use the Gelfand formula $\rho (P_k) =\displaystyle\lim_{r \to \infty}\| P_k^r \|^{1/r}$ for the spectral radius of a matrix and follow the ideas of \cite{nonBT}. Let $C_1$ and $C_2$ be two cycles in $G$ of lengths $\ell_1$ and $\ell_2$, respectively, with at least one different vertex and with $\ell_1 \geq \ell_2 > k$. Since $G$ is strongly connected, both cycles are connected by two paths $p_1$ and $p_2$ of lengths $d_{1}$ and $d_{2}$, respectively, that connect one vertex in $C_1$ with a vertex in $C_2$ in both orientations. 
	
\begin{figure}[h]
	\begin{center}
		\begin{tikzpicture}[scale=0.5]
		\node[] (1) at (1,0) {};
		\draw (1,0) node {$\bullet$};
		\node[] (2) at (6,0) {};
		\draw (6,0) node {$\bullet$};
		
		\draw[dashed,->] (1) to[bend left] (2) ;
		\draw (3.5,1.5)  node {$p_1$} ;
		\draw[dashed,->] (2) to[bend left] (1) ;
		\draw (3.5,-1.5)  node {$p_2$} ;
		\draw (1,0) arc (0:359:2)  ;
		\draw (-1,0) node {$C_1$};
		\draw (6,0) arc (359:0:-1.5) ;
		\draw (7.5,0) node {$C_2$};
		\end{tikzpicture}
	\end{center}
	\vspace{-0.5cm}
	\caption{Scheme for the proof of Theorem \ref{thm:directedPk} }
\end{figure}
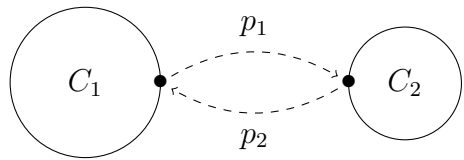

{Given an} integer $t$ such that $t\geq 2(\ell_1 + \ell_2 + \frac{d_1 + d_2}{2} )$ {there exists an unique integer $s$ satisfying}
	$$2s\left(\ell_1 + \ell_2 + \frac{d_1 + d_2}{2} \right) \leq t < (2s+1)\left(\ell_1 + \ell_2 + \frac{d_1 + d_2}{2} \right) . $$
	We will show that there is at least one entry of $P_k^t$ that is bounded below by $2^{2s}$ so that then $$\rho (P_k) \geq \displaystyle\lim_{t \to \infty}(2^{2s})^{1/t}=2^{1/ \left( \ell_1 + \ell_2 + \frac{d_1 + d_2}{2}\right)}>1.$$

	First, we consider the case in which $C_1$ and $C_2$ have one vertex in common and we can consider $d_1=d_2=0$. We note that one of the entries of $P_k^t$ is the number of non-$k$-cycling walks of length $t$ starting and ending somewhere in $C_1$. Moreover, the number of those non-$k$-cycling walks that go precisely $2s$ times
	around $C_1$ and $2s$ times around $C_2$ is $\mathcal{PR}_{2s,2s}^{4s}$\footnote{$\mathcal{PR}_{n_1,n_2,\ldots,n_m}^{n}=\frac{n!}{n_1!n_2!\cdots n_m!}$ is the number of permutations of $n$ objects of which $n_1$ are of one kind, $n_2$ are of second kind,$\ldots$, and $n_m$ are of $m^{\text{th}}$ kind.} and $$\mathcal{PR}_{2s,2s}^{4s}=\frac{(4s)!}{(2s)!(2s)!}=\frac{4s}{2s} \cdot\frac{4s-1}{2s-1}\, \cdots \,\frac{2s+1}{1}\geq 2^{2s}.$$

We now assume that  $C_1$ and $C_2$ do not have one vertex in common and consider again non-$k$-cycling walks of length $t$ starting and ending somewhere in $C_1$. Let us count the number of those non-$k$-cycling walks that go precisely $2s$ times
around $C_1$, $2s$ times around $C_2$, $s$ times through $p_1$ and $s$ times through $p_2$.
{Note that, in our construction, it is possible that $p_1$ and $p_2$ can potentially share vertices with $C_1$ and $C_2$. However, this poses no problem with respect to backtracking. On the other hand, note that the paths $p_1$ and $p_2$ could also possibly share the same vertices between themselves. In this latter case,} in order to avoid backtracking, we will go at least one time around $C_2$ each time we go through $p_1$ and at least one time around $C_1$ each time we go through $p_2$.

 That is, we are fixing $s$ of the times we go around $C_2$ and $s-1$ of the times around $C_1$. Then, for the others $2s-s=s$ times around $C_2$ we have $\mathcal{CR}_{s}^{s}={2s-1 \choose s}$ different options, and for the others $2s-(s-1)=s+1$ times around $C_1$ we have $\mathcal{CR}_{s+1}^{s+1}={2s+1 \choose s+1}$ different options \footnote{$\mathcal{CR}_{m}^{n}={n+m-1 \choose n}$ is the number of combinations of $m$ objects, taken $n$ at a time, when repetition of objects is allowed.}. We note that 
$$\mathcal{CR}_{s}^{s}=\frac{(2s-1)!}{(s-1)!s!}=\frac{2s-1}{s-1} \cdot\frac{2s-2}{s-2}\, \cdots \,\frac{s+1}{1}\geq 2^{s-1}$$ and
$$\mathcal{CR}_{s+1}^{s+1}=\frac{(2s+1)!}{(s+1)!s!}=\frac{2s+1}{s+1} \cdot\frac{2s}{s}\, \cdots \,\frac{s+1}{1}\geq 2^{s+1}.$$ Therefore, 
$$ \mathcal{CR}_{s}^{s}\,\mathcal{CR}_{s+1}^{s+1}\geq 2^{s-1} 2^{s+1} = 2^{2s}.$$
\end{proof}

\subsection{Radius of convergence}

In this subsection we aim to state and prove Theorem \ref{th:radius}, where we characterize the radius of convergence of a digraph $G$ in terms of the number of cycles in the undirectizations of its strongly connected components. Our first step is to study the eigenvalues of the directed deformed graph Laplacian.

\begin{theorem}\label{th:eigenvalues}
Let $G$ be a strongly connected digraph, and consider $H$ the undirectization of $G$. Let $M(t)$ be the directed deformed graph Laplacian of $G$. Then the following statements hold:
\begin{itemize}
\item[$\rm (a)$] If $H$ is a tree then $G$ is an undirected tree and the only finite eigenvalues of $M(t)$ are $\pm 1$.
\item[$\rm (b)$] If $H$ has precisely one cycle of length $\ell$ then $G$ does not have any two directed cycles with at least one different vertex and the only finite eigenvalues of $M(t)$ are the $\ell$-th roots of unity.
\item[$\rm (c)$] If $H$ has at least two cycles then $G$ has at least two directed cycles with at least one different vertex and $M(t)$ has at least one real positive eigenvalue in the open interval $(0,1)$.
\end{itemize}
\end{theorem}

\begin{proof} Note that $H$ is connected since $G$ is strongly connected.
\begin{itemize}
\item[$\rm (a)$] Since $H$ is a tree, $G$ is an undirected tree  by Lemma \ref{lem:undirectedtree}. It is known \cite{GHN} that the deformed graph Laplacian of undirected trees only have $\pm 1$ as finite eigenvalues.
\item[$\rm (b)$] Since $H$ has precisely one cycle, $H$ is either a cycle graph or is a cycle plus some tails. In the latter case, then $G$ must have the same \emph{undirected} tails, as otherwise we contradict strong connectedness, as in the proof of Lemma \ref{lem:undirectedtree}. On the other hand, by Theorem \ref{thm:leaves} (recursively applied), we can without loss of generality remove all the tails. Thus, we can assume that $H$ is a cycle graph with $\ell$ vertices. $G$ is a subgraph of $H$, but since it must contain at least one cycle this implies that it is also a supergraph of one of the two directed cycle graphs with $\ell$ vertices. Then, it is clear that the graph of the Hashimoto matrix $B$ of $G$ is either the union of two directed cycles (if $G=H$) of length $\ell$ or the union of one directed cycle of length $\ell$ and possibly some directed path graphs. Therefore, the eigenvalues of $B$ are either $0$ or $\ell$-th roots of unity. Since the eigenvalues of $M(t)$ are either $\pm 1$ or the inverses of the eigenvalues of $B$ by Ihara's theorem for digraphs (see \cite{tarfulea} or Corollary \ref{cor:Ihara}), we deduce that the finite eigenvalues of $M(t)$ are all unimodular, and in particular, there is no eigenvalue on $(0,1)$.
\item[$\rm (c)$]
Since $H$ has at least two cycles then, by Lemma \ref{lem:HimpliesG}, $G$ has at least two directed cycles with at least one different vertex. Then, by Theorem \ref{thm:directedPk}, $\rho(B)>1$, where $B$ is the Hashimoto matrix.
By the Perron-Frobenius theorem, $\rho(B)$ is an eigenvalue of $B$. In addition, we know that a complex number not equal to $\pm 1$ is a finite eigenvalue of $M(t)$ if and only if its inverse is an eigenvalue of $B$ by Ihara's theorem for digraphs (stated in Corollary \ref{cor:Ihara}). Therefore, $1/\rho(B)<1$ is an eigenvalue of $M(t)$.
\end{itemize}
\end{proof}

\begin{remark}\rm
    Note that the extreme case of $G$ being just an isolated node corresponds to $A=S=D=0$ and $M(t)=1-t^2$.
\end{remark}
 In summary, we have therefore proved the following result.

\begin{theorem} Let $M(t)$ be the directed deformed graph Laplacian of a digraph $G$. Then $M(t)$ has at least one real positive eigenvalue inside the unit open circle if and only if $G$ contains at least one strongly connected component whose undirectization has at least two cycles. Otherwise, the eigenvalues of $M(t)$ are all unimodular.
\end{theorem}

By using Theorem \ref{th:eigenvalues}, the following Theorem \ref{th:radius} gives the complete characterization of the radius of convergence of the generating function $$\sum_{k=0}^\infty t^k p_k(A)=(1-t^2)M(t)^{-1}.$$

\begin{theorem}\label{th:radius}
Let $M(t)$ be the {directed} deformed graph Laplacian of a digraph $G$. Let $\mu$ with $0 < \mu \leq 1$ be the smallest positive eigenvalue of $M(t)$, and let $r$ be the radius of convergence of $\sum_{k=0}^\infty t^k p_k(A)=(1-t^2)M(t)^{-1}$. Assume that $G$ has $s$ strongly connected {or single node} components and denote by $H_1, \dots, H_s$ the undirectizations of each of them. Then the following statements hold:
\begin{itemize}
\item[$\rm (a)$] If all the undirectizations $H_i$ are trees {or isolated nodes,} then $1 = \mu < r = \infty$. 
\item[$\rm (b)$] If no one of the undirectizations $H_i$ has more than one cycle but at least one of them has precisely one cycle then $1 = \mu = r$.
\item[$\rm (c)$] If at least one undirectization $H_i$, for some $1 \leq i \leq s$, has more than one cycle then $1 > \mu = r$.
\end{itemize}
\end{theorem}

\begin{proof}
	\begin{itemize}
		\item[$\rm (a)$] By Theorem \ref{th:eigenvalues}$\rm (a)$, if all the undirectizations $H_i$ are trees {or isolated nodes,} then the only finite eigenvalues of all the {directed} deformed graph Laplacians $M_{i}(t)$ of the strongly connected components of $G$ are $\pm 1$. Then, { by Proposition \ref{prop:algebraic}}, the only finite eigenvalues of $M(t)$ are $\pm 1$. Since each $H_i$ has no cycles, by Corollary \ref{cor:defective}, $ 1$ is a semisimple eigenvalue of each $M_{i}(t)$. Therefore, by Remarks \ref{rem:algebraic} and \ref{rem:geometric1}, $1$ is also a semisimple eigenvalue of $M(t)$. Thus, $1$ is not a pole of the generating function. Then, we see by contradiction that $-1$ is not a pole of the generating function either since, in such a case, $r=1$ and by Abel's Theorem on power series \cite{Abel} this implies that $1$ is a pole. Therefore, $r=\infty$.
		
		\item[$\rm (b)$] By Theorem \ref{th:eigenvalues}$\rm (b)$, if no one of the undirectizations has more than one cycle but at least one of them $H_i$ has precisely one cycle of length $\ell$ then the only finite eigenvalues of the deformed graph Laplacian $M_{i}(t)$ associated to the corresponding strongly connected component are the $\ell$-th roots of unity. Then, { by Proposition \ref{prop:algebraic}}, the eigenvalues of $M(t)$ are the $\ell$-th roots of unity. In addition, by Corollary \ref{cor:defective}, $1$ is a defective eigenvalue of $M_{i}(t)$. Therefore, $1$ is also a defective eigenvalue of $M(t)$ and, thus, $1$ is a pole of the generating function.
		
		\item[$\rm (c)$] By Theorem \ref{th:eigenvalues}$\rm (c)$, if at least one undirectization $H_i$ has more than one cycle then the {directed} deformed graph Laplacian $M_{i}(t)$ associated to the corresponding strongly connected component {of $G$} has at least one real positive eigenvalue $\lambda$ in the open interval $(0,1)$. Therefore, $\lambda$ is also an eigenvalue of $M(t)$ { by Proposition \ref{prop:algebraic}}. Since $\lambda\neq 1$, we conclude that $\lambda$ is also a pole of the generating function and the result follows.
	\end{itemize}
\end{proof}

\section{BTDW centrality}\label{sec:BTDW}

In this section we consider walks on an unweighted digraph with no loops where backtracking is attenuated rather than eliminated, which are called Backtrack-Downweighted Walks (BTDWs) \cite{GHN}. That is, if a walk of length $k$ backtracks $m$ times, it is weighted by $t^k \omega ^m$, where $0\leq \omega\leq 1$ is a parameter. Note that $\omega=0$ corresponds to NBTW while $\omega=1$ to classical Katz {\cite{Katz}}. Denote by $q_k(A,\omega)$ the matrix whose entry $(q_k(A,\omega))_{ij}$ counts the number of BTDWs. It was proved in \cite{inprep} that the following expression for the generating function holds:

\begin{equation}\label{eq:genfunbtdw}
    \sum_{k=0}^{\infty} t^k q_{k}(A,\omega)=(1-\tau^2 t^2)(I - A t + (D- \tau I )\tau t^2 + (A-S) \tau^2t^3)^{-1},
\end{equation}
where $\tau:= 1-\omega$ and $D,$ $A$ and $S$ are as in Section \ref{sec:NBTW}. If the corresponding graph is undirected then the generating function is simplified as $S=A$. We define the polynomial matrix
\begin{equation*}\label{eq:BTDWdeformedgraphlaplacian}
M_{\tau}(t):=I - A t + (D- \tau I )\tau t^2 + (A-S) \tau^2t^3,
\end{equation*}
which we call the {\em $\tau$--deformed graph Laplacian} of $G$, and we will study the radius of convergence of the generating function and the spectrum of $M_{\tau}(t)$ with respect to the underlying digraph. Note that $\tau=1$ is the NBTW case.

\subsection{Eigenvalues of the $\tau$--deformed graph Laplacian: Algebraic, geometric and partial multiplicities}

Let $G$ be a digraph and let $M_{\tau}(t)$ be its $\tau$--deformed graph Laplacian. Assume that the adjacency matrix $A$ of $G$ can be written in the form  
\[ A= \begin{bmatrix}
A_1 & E\\
0 & A_2
\end{bmatrix} .\] Then the $\tau$--deformed graph Laplacian of $G$ can be written as
\[  M_{\tau}(t) = \begin{bmatrix}
M_{\tau}^{1}(t) & (\tau^2 t^3-t) E\\
0 & M_{\tau}^{2}(t) 
\end{bmatrix} .\]

Therefore, as in the NBTW case, we note that in order to studying the algebraic multiplicities of eigenvalues of the $\tau$--deformed graph Laplacian $M_{\tau}(t)$, there is no loss of generality in assuming that the graph $G$ is strongly connected.

For the eigenvalues $\pm \frac{1}{\tau}$, the geometric multiplicities is also the sum of the geometric multiplicities for each $\tau$--deformed graph Laplacians $M_{\tau}^{i}(t)$ of each strongly connected component since $M(\pm \frac{1}{\tau})= M_{\tau}^1(\pm \frac{1}{\tau}) \oplus M_{\tau}^2(\pm \frac{1}{\tau})$. In addition, we have the following result.

\begin{proposition}
	Let $M_{\tau}(t)$ be the $\tau$--deformed graph Laplacian of a digraph $G$. Then, $\frac{1}{\tau}$ is always an eigenvalue, and its geometric multiplicity is the number of connected components of $G_U$, the undirected part of $G$. Moreover, $-\frac{1}{\tau}$ is an eigenvalue if and only if $G_U$ has at least one bipartite connected component. The geometric multiplicity of $-\frac{1}{\tau}$ is equal to the number of such components.
\end{proposition}

\begin{proof} Denoting by $M_{\tau}(t)$ and $M_{\tau}^{U}(t)$ the $\tau$--deformed graph Laplacians of $G$ and $G_U$, respectively, we have that
	\begin{equation}\label{eq:2}
	M_{\tau}(t) = M_{\tau}^U(t) + (\tau^2t^3-t) (A-S).
	\end{equation} Let us now denote by $L:=D-S$ and $Q:=D+S$, respectively, the graph Laplacian and the signless Laplacian of $G_U$, the undirected part of $G$. Then, by \eqref{eq:2}, $M_{\tau}(\frac{1}{\tau}) = M_{\tau}^{U}(-\frac{1}{\tau}) = \frac{1}{\tau}L$. Let $e$ be the vector of all ones. Then $Le=0$ and the nullity of $L$, i.e., the geometric multiplicity of $\frac{1}{\tau}$, is the number of connected component of $G_U$. Moreover, we have that $M_{\tau}(-\frac{1}{\tau})=M_{\tau}^{U}(-\frac{1}{\tau})=\frac{1}{\tau}Q$ by \eqref{eq:2}, and $Q$ is singular if and only if $G_U$ has at least one bipartite component. In addition, the geometric multiplicity is the number of such bipartite components. 	
\end{proof}

\begin{proposition}\label{prop:defective2} Let $G$ be a digraph and let $M_{\tau}(t)$ be its $\tau$--deformed graph Laplacian. Let $d$ and $d_U$ be defined as in \eqref{def:total} and \eqref{def:totalrecip}. Then $\frac{1}{\tau}$ is a defective eigenvalue of $M_{\tau}(t)$ if and only if $2 d - d_U = 2\tau n$.
\end{proposition}

\begin{proof}
	$\frac{1}{\tau}$ is a defective eigenvalue of $M_{\tau}(\lambda)$ if and only if there exists a vector $v$ such that $M_{\tau}(\frac{1}{\tau})v+M_{\tau}'(\frac{1}{\tau})e=0$. Since $M_{\tau}'(1)=2L+ 2A - 2\tau I - S$ then $M_{\tau}(\frac{1}{\tau})v+M_{\tau}'(\frac{1}{\tau})e=0$ for some $v$ if and only if  $0 = e^T (2A - 2\tau I - S) e  = 2d -2\tau n - d_U$.
\end{proof}

For a single node, the $\tau$-deformed graph Laplacian is simply $1-\tau^2 t^2$, which evidently has always eigenvalues $\pm 1/\tau$. By Remark \ref{rem:average} and Proposition \ref{prop:defective2}, we obtain the following corollary on strongly connected components.

\begin{corollary} Let $G$ be a strongly connected digraph and let $M_{\tau}(t)$ be its $\tau$--deformed graph Laplacian. Let $H$ be the undirectization of $G$. Then
	\begin{itemize}
		\item[\rm (a)] If $\tau>1$, $1/\tau$ is a defective eigenvalue of $M_{\tau}(t)$ if and only if $H$ has at least two cycles;
		\item[\rm (b)] If $\tau =1 $, $1/\tau$ is a defective eigenvalue of $M_{\tau}(t)$ if and only if $H$ has precisely one cycle;
		\item[\rm (c)] If $\tau<1$, $1/\tau$ is a defective eigenvalue of $M_{\tau}(t)$ if and only if $H$ is a tree.
	\end{itemize}
	
\end{corollary}

\begin{proposition} Let $G$ be a digraph and let $M_{\tau}(t)$ be its $\tau$--deformed graph Laplacian.
	If $\lambda \neq \pm \frac{1}{\tau}$ is a finite eigenvalue of only one $\tau$--deformed graph Laplacian $M_{\tau}^{i}(t)$ associated with a strongly connected component $G_i$ of $G$, then its geometric multiplicity as an eigenvalue of $M_{\tau}(t)$ and as an eigenvalue of $M_{\tau}^{i}(t)$ is the same.
\end{proposition}
\begin{proof}
	If $\lambda \neq \pm \frac{1}{\tau}$ is an eigenvalue only for $M_{\tau}^1(t)$ then $M_{\tau}^1(\lambda) U=0$ if and only if
	\[  \begin{bmatrix}
	M_{\tau}^1(\lambda) & (\tau^2\lambda^3-\lambda) E\\
	0 & M_{\tau}^2(\lambda)
	\end{bmatrix} \begin{bmatrix}
	U\\
	0
	\end{bmatrix} = 0\]
	for any full column rank matrix $U$, which columns would be eigenvectors associated with $\lambda$. On the other hand, if $\lambda$ is an eigenvalue only for $M_{\tau}^2(t)$ then $M_{\tau}^2(\lambda)V=0$ if and only if
	\[  \begin{bmatrix}
	M_{\tau}^1(\lambda) & (\tau^2\lambda^3-\lambda) E\\
	0 & M_2(\lambda)
	\end{bmatrix} \begin{bmatrix}
	(\lambda- \tau^2\lambda^3) M_{\tau}^1(\lambda)^{-1} E V\\
	V
	\end{bmatrix} = 0,\]
	for any full column rank matrix $V$.
\end{proof}

However, if $\lambda \neq \pm \frac{1}{\tau}$ is a finite eigenvalue of both $\tau$--deformed graph Laplacians $M_{\tau}^1(t)$ and $M_{\tau}^2(t)$ of each strongly connected component, then its geometric multiplicity as an eigenvalue of $M_{\tau}(t)$ may not be equal to the sum of the geometric multiplicities (recall Example \ref{ex_multiplicities}).

\begin{remark}\rm
In contrast to the NBTW case (i.e., $\tau=1$), pruning reciprocal leaves does not preserve finite eigenvalues of the $\tau$--deformed graph Laplacian if $\tau\neq 1$. Let us see a simple example. 
\end{remark}

\begin{example}
    Consider the following digraph
    
	
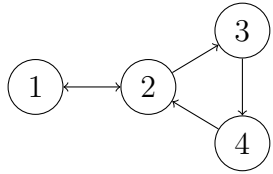
\begin{figure}[H]
		\begin{center}
			\begin{tikzpicture}[scale=0.5]
			\node[draw,circle] (1) at (1,0) {$1$};
			\node[draw,circle] (2) at (4,0) {$2$};
			\node[draw,circle] (3) at (6.5,1.5) {$3$};
			\node[draw,circle] (4) at (6.5,-1.5) {$4$};
			\draw[<->] (1) -- (2);
			\draw[->] (2) -- (3);
			\draw[->] (3) -- (4);
			\draw[->] (4) -- (2);
		\end{tikzpicture}
		\end{center}
		\vspace{-0.5cm}
		\end{figure}
  Denote by $M_{\tau}(t)$ its $\tau$--deformed graph Laplacian and let $\widetilde M_{\tau}(t)$ be the $\tau$--deformed graph Laplacian obtained by removing the reciprocal leaf $1 \leftrightarrow 2$. Then, for $\tau=1/2$, it can be proved that the Smith form of $M_{\tau}(t)$ is $\diag(1, t^2 - 4 , t^2 - 4, t^7 + t^4 - 16t^3 - 8t^2 + 16  )$. However, the Smith form of $\widetilde M_{\tau}(t)$ is $\diag(t^2 - 4,t^2 - 4, t^5 - 4t^3 - t^2 + 4)$, and the last invariant polynomials $t^7 + t^4 - 16t^3 - 8t^2 + 16$ and $t^5 - 4t^3 - t^2 + 4$ have disjoint sets of roots.
\end{example}

\subsection{Ihara's theorem for digraphs including downweighting}

In this section, we extend Ihara's Theorem \cite{ihara} to the case of backtrack-downweighted walks for unweighted digraphs. {We first introduce the source and target matrices associated with a digraph.

\begin{definition}\label{def:source_target} Let $G=(V,E)$ be a digraph (weighted or unweighted) of order $n$ with $m$ edges. The source matrix $L\in \{0,1\}^{m \times n}$ of $G$ is defined as:
	$$L_{ej}= \left\{ \begin{array}{lcc}
	1 &   \text{ if the edge }  e \text{ has the form } e=(j,\star)    \\
	0& otherwise 
	\end{array}
	\right.,$$
	and the target matrix $R\in \{0,1\}^{m \times n}$ of $G$ is defined as:
	$$R_{ej}= \left\{ \begin{array}{lcc}
	1 &   \text{ if the edge }  e \text{ has the form } e=(\star,j)    \\
	0& otherwise 
	\end{array}
	\right. ,$$
	for $e=1,\ldots,m$ and $j=1,\ldots,n.$
\end{definition}} The source and target matrices give a factorization for the adjacency matrix $A$ of {an unweighted digraph} $G$, as it can be proved that $A=L^TR$. The matrix $W:=RL^T$ is called the line graph adjacency matrix. We also set 
	\begin{equation}\label{eq:matrices}
\Delta:=W \circ W^T, \quad \Omega:=\Delta^2 \quad\text{and} \quad B:=W-\Delta.
	\end{equation}
	The matrix $B$ thus defined is the Hashimoto matrix \cite{Hashimoto} and coincides also with the non-$2$-cylcing matrix $P_2$ in Subsection \ref{sub:nonkcycling} and in \cite{nonBT}: see, e.g., \cite[Theorem 5.3]{nonBT}. Before stating the main result in this section (Theorem \ref{th:det}), we need the following lemmas.

	\begin{lemma}
		The matrix $\Omega$ in \eqref{eq:matrices} is diagonal and 
			$$\Omega_{ee}= \left\{ \begin{array}{lcc}
		1 &   \text{ if } e \text{ is reciprocated}   \\
		0& otherwise 	\end{array}
		\right. .$$
		 In particular, $\Omega_{ee}=1=\Omega_{ff}$ if and only if $ \Delta_{ef}=1=\Delta_{fe}$.
	\end{lemma}
	\begin{proof}
		We have $\Delta_{ef}=1$ if and only if there is a walk of length $2$ $e \rightarrow f$ and there is a walk of length $2$ $f \rightarrow e$, i.e., if and only if $f$ is the reciprocal edge of $e$. Moreover, $\Delta_{ef}=0$ in all other cases. Hence, $\Omega_{ee} = \sum_f \Delta_{ef} \Delta_{fe}=1$ if $e$ is reciprocated and $\Omega_{ee}=0$ otherwise. 
	\end{proof}
	
	\begin{lemma}\label{lem:powers}Consider the matrix $\Delta$ in \eqref{eq:matrices}. For all $k\geq1$, 
		$$ \Delta^k = \begin{cases}
		\Delta \ \mathrm{if} \ k \ \mathrm{odd},\\
		\Omega \ \mathrm{if} \ k \ \mathrm{even}.
		\end{cases}$$
	\end{lemma}
	\begin{proof}
		It suffices to show $\Delta^3=\Delta$, then the other cases follow by induction on $k$. To this goal note that $\Delta^3 = \Delta \Omega$. Indeed, $(\Delta^3)_{ef} = \Delta_{ef} \Omega_{ff}=\Delta_{ef}$.
	\end{proof}
	
	\begin{lemma}\label{lem:D_S} Let $L$ and $R$ be the source and target matrices of a digraph $G$, respectively. Consider the matrices $\Delta$ and $\Omega$ in \eqref{eq:matrices} and the matrices $S$ and $D$ defined in Section \ref{sec:preliminary}. Then
		$$L^T \Delta R = D\quad \text{and} \quad L^T \Omega R = S.$$
	\end{lemma}
	\begin{proof}First, we have that
		$$L^T \Delta R = L^T W R - L^T B R = A^2 - p_2(A) = D.$$
		On the other hand $$(L^T \Omega R)_{ij} = \sum_{k} L_{ki} \Omega_{kk} R_{kj} = \begin{cases} 1 \ \mathrm{if}\ ij\in E \ \mathrm{and}\ ji \in E\\
		0 \ \mathrm{otherwise}
		\end{cases},$$
		which by definition is also equal to $S_{ij}$.
	\end{proof}
	
	\begin{lemma}\label{lem:det}{Let $G$ be a digraph and let $\Delta$ be the matrix} in \eqref{eq:matrices}, then
		$$\det(t I - \Delta) = t^a (t^2-1)^b$$ where $a$ is the number of unreciprocated edges and $b$ is the number of (undirected) reciprocal edges { in $G$}.
	\end{lemma}
	\begin{proof}
		By the previous lemma it is clear that the eigenvalues of $\Delta$ must belong to $\{0,1,-1\}$. It is also clear that the multiplicity of $0$ is the number of unreciprocated edge, as that is also equal to the number of zero diagonal elements of $\Omega$. Finally, up to a permutation of the vertices so that they are ordered mentioning first reciprocal edges, $$\Delta = \underbrace{ \begin{bmatrix}
			0 & 1\\
			1 & 0
			\end{bmatrix} \oplus \dots \oplus \begin{bmatrix}
			0 & 1\\
			1 & 0
			\end{bmatrix}}_{b \ \mathrm{copies}} \oplus \underbrace{0 \oplus \dots \oplus 0}_{a \ \mathrm{copies}}.$$
	\end{proof}

{
\begin{lemma}\label{lem:downweighting} Let $G$ be a digraph and let $M_{\tau}(t)$ be its $\tau$--deformed graph Laplacian. Consider $L$ and $R$ the source and target matrices of $G$, respectively, and the matrix $\Delta$ in \eqref{eq:matrices}. Then,

$$I - t L^T[I+\tau t \Delta]^{-1} R= \dfrac{1}{1-\tau^2 t^2}\,M_{\tau}(t).$$

\end{lemma}

\begin{proof} 
		By Lemma \ref{lem:powers}, we have that $$ (I+t \tau \Delta) \left( I + \Omega \dfrac{t^2 \tau^2}{1-t^2 \tau^2} - \Delta \dfrac{t \tau}{1-t^2 \tau^2}\right) = I.$$ Then, taking into account Lemma \ref{lem:D_S} and the fact that $A=L^TR$, we obtain $$ L^T(I+\tau t \Delta)^{-1}R = A + S \frac{\tau^2 t^2}{1-\tau^2 t^2} - D \frac{\tau t}{1-\tau^2 t^2}.$$
  Therefore,
  $$I - t L^T[I+\tau t \Delta]^{-1} R=\dfrac{1}{1-\tau^2 t^2} [ I-tA+t^2\tau( D- \tau I)+\tau^2 t^3(A-S)]= \dfrac{1}{1-\tau^2 t^2}\,M_{\tau}(t).$$ 
  
\end{proof}
}

	\begin{theorem}\label{th:det} Let $G$ be a digraph of order $n$ and let $M_{\tau}(t)$ be its $\tau$--deformed graph Laplacian. Assume that $G$ has $b$ (undirected) reciprocal edges. Then
		\[\det [I - t  (\tau B + (1-\tau) W) ]  =  (1-\tau^2 t^2)^{b-n} \det M_{\tau}(t),  \]
		for all $0 \leq \tau \leq 1.$
	\end{theorem}
	\begin{proof}
		Since $B=R L^T-\Delta$, then $\tau B + (1-\tau) W = W - \tau \Delta$. {Taking into account that $W=RL^{T}$,} by the generalized matrix determinant lemma, we have that
		\[ \det (I - t W + t \tau \Delta ) = \det (I + t \tau \Delta) \det ( I - t L^T [I + t \tau \Delta]^{-1} R ).   \]
		{Finally, the result follows by noting that, on one hand, by Lemma \ref{lem:det},}
		\[ \det (I + \tau t \Delta) = (-\tau t)^{a+2b} \det( {-\frac{1}{\tau t}}I - \Delta) = (-\tau t)^{a+2b} (-\tau t)^{-a} ((\tau t)^{-2}-1)^b = (1-\tau^2 t^2)^b ;\]
		{and, on the other hand, by Lemma \ref{lem:downweighting},}
				\[ \det ( I - t L^T[I+\tau t \Delta]^{-1} R) = (1-\tau^2 t^2)^{-n} {\det M_{\tau}(t)}. \]
	\end{proof}

The following corollary follows from Theorem \ref{th:det}.

\begin{corollary}\label{cor:inv} Let $G$ be a digraph and let $M_{\tau}(t)$ be its $\tau$--deformed graph Laplacian. Denote $V(\tau):=\tau B + (1-\tau) W$ and let $\lambda \neq \pm 1/\tau$. Then $\lambda$ is an eigenvalue of $M_{\tau}(t)$ if and only if $1/\lambda$ is an eigenvalue of $V(\tau)$.
\end{corollary}
	
	Note that, as further corollaries of Thereom \ref{th:det} for $\tau=1$ and $\tau=0$ respectively,  we obtain two well-known results: Ihara's Theorem {for digraphs (that can be found in \cite{tarfulea} but stated in a different form involving zeta functions)}, and a particular case of Flanders' Theorem \cite{flanders}.
	
	\begin{corollary}[Ihara's Theorem {for digraphs}]\label{cor:Ihara} Let $G$ be a digraph of order $n$ with associated {directed} deformed graph Laplacian $M(t)$ and Hashimoto matrix $B$. Assume that $G$ has $b$ (undirected) reciprocal edges. Then
		\[ \det (I - t B) = (1-t^2)^{b-n} \det M(t).\]
	\end{corollary}

	\begin{corollary}[A special case of Flanders' Theorem] Let $G$ be a digraph with adjacency matrix $A$ and line graph adjacency matrix $W$. Then
	$ \det (I - t   W )  =   \det (I - t A). $
\end{corollary}

\subsection{Radius of convergence}
In this subsection, we study the radius of convergence of \eqref{eq:genfunbtdw}.
\begin{theorem}\label{th:eigenvalues2}
	Let $G$ be a strongly connected digraph, and consider $H$ the undirectization of $G$. Let $M_{\tau}(t)$ be the $\tau$--deformed graph Laplacian of $G$ with $0<\tau\leq 1$. If $H$ has at least two cycles then $M_{\tau}(t)$ has at least one real positive eigenvalue in the open interval $(0,1/\tau)$.
\end{theorem}

\begin{proof} Consider the matrix $V(\tau)$ in Corollary \ref{cor:inv} and note that $\rho(V(\tau))\geq \rho(B)$ if $0<\tau\leq 1$. In addition, by Theorem \ref{thm:directedPk}, if $H$ has at least two cycles then $\rho(B)>1\geq \tau$, where $B$ is the Hashimoto matrix. Therefore,
	\begin{equation}\label{eq:ineq}
\rho(V(\tau))\geq \rho(B)>1\geq \tau.
	\end{equation}
	By the Perron-Frobenius theorem \cite[Theorem 6.2.1]{Friedland}, $\rho(V(\tau))$ is an eigenvalue of $V(\tau)$, and, thus, by Corollary \ref{cor:inv}, $1/\rho(V(\tau))$ is an eigenvalue of $M_{\tau}(t).$ Taking into account \eqref{eq:ineq},  $1/\rho(V(\tau))<1/\tau$.
\end{proof}

\begin{theorem} Let $M_{\tau}(t)$ be the $\tau$--deformed graph Laplacian of a digraph $G$. Then $M_{\tau}(t)$ has at least one real positive eigenvalue inside the open circle of radius $1/\tau$ if $G$ contains at least one strongly connected component whose undirectization has at least two cycles.
\end{theorem}

We do not have a complete characterization of the radius of convergence for the BTDW case, but the following result that establishes how the radius of convergence is in a very common case in application.

\begin{theorem}\label{th:radius2}
Let $M_{\tau}(t)$ be the $\tau$--deformed graph Laplacian of a digraph $G$. Let $\mu$ with $0 < \mu \leq 1/\tau$ be the smallest positive eigenvalue of $M_{\tau}(t)$, and let $r$ be the radius of convergence of $\sum_{k=0}^\infty t^k q_k(A,\omega)=(1-\tau^2 t^2)M_\tau(t)^{-1}$ as in \eqref{eq:genfunbtdw}. Assume that $G$ has $s$ strongly connected {or single node} components and denote by $H_1, \dots, H_s$ the undirectizations of each of them. Then if at least one undirectization $H_i$, for some $1 \leq i \leq s$, has more than one cycle then $ \mu = r<1/\tau$.

\end{theorem}

\section{The radius of convergence and Ihara's Theorem in the case of weighted {di}graphs}\label{sec:weighted}

In this section, we study {the radius of convergence of the generating function for the case of} weighted digraphs.

 Let $G=(V,E)$, with $\#V=:n$ and $\#E=:m$, be a weighted directed graph with adjacency matrix $A$. In this setting, each edge is associated with a positive weight; recall from Section \ref{sec:preliminary} that the weight a walk of length $r$, say, $(e_1,\dots,e_r)$, is defined as the product $\prod_{i=1}^r \omega(e_i)$ of the weights of the $r$ edges (counted with repetitions) that are used in the walk. The definition of the matrices $p_k(A)$ can be accordingly generalized so that the $(i,j)$ element of $p_k(A)$ is equal to the weighted sum of all non-backtracking walks of length $k$, where the weight of each walk is defined as above. Note that the unweighted case can be seen as the limit where the weight of all edges (and hence all walks) are equal to $1$.

In \cite[Corollary 4.7]{ryan} it was shown that, in this more general context, the generating function \eqref{eq:genfun} converges to
\begin{equation}\label{eq:genfunweighted}
    \Phi(A,t) = \sum_{k=0}^\infty t^k p_k(A) = I + t L^T \sqrt{Z} (I-tV)^{-1} \sqrt{Z} R,
\end{equation}
{where the matrices $L$ and $R$ are the source and target matrices, respectively, in Definition \ref{def:source_target}; and the matrices $Z$ and $V$ are constructed as in Definition \ref{def:new} below \cite[Section 4]{ryan}}. To state Definition \ref{def:new}, we denote by $\R_+:=]0,+\infty[$ the set of positive real numbers and by $\R_\geq:=[0,+\infty[$ the set of non-negative real numbers.
\begin{definition}\label{def:new}  
    \begin{itemize}
    \item[$\rm(a)$] $Z \in \R_\geq^{m \times m}$ is a diagonal matrix such that $Z_{ee} \in \R_+$ is the weight of the directed edge $e$. Note that $Z$ is positive definite and that by construction we have $A=L^TZR$.
    \item[$\rm(b)$]  Let $B \in \R_\geq^{m \times m}$ be the weighted version of the Hashimoto matrix \cite{Hashimoto}. In particular $B_{ef}>0$ if and only if $(e,f)$ is a NBTW of length $2$ and, in this case, $B_{ef}=\omega(e)\omega(f)$ is the weight of such walk. Then we set
    \begin{equation}\label{eq:Hashimoto_weight}
        V:=B^{\circ 1/2}
    \end{equation}
\end{itemize}
\end{definition}

Let $r$ be the radius of convergence of \eqref{eq:genfunweighted}. Via \eqref{eq:genfunweighted}, it was proved in \cite[Corollary 4.8]{ryan} that $\rho(V)^{-1} \leq r$. The goal of this section is to improve on \cite[Corollary 4.8]{ryan} and show that in fact $\rho(V)^{-1}$ is not only a lower bound but actually equal to the radius of convergence.

\begin{theorem}\label{thm:weightedradius}
     Let $G$ be a weighted digraph, and let $V$ be the matrix in \eqref{eq:Hashimoto_weight}. Then the radius of convergence $r$ of the generating function in \eqref{eq:genfunweighted} is $r=\rho(V)^{-1}$. Here, we allow $\rho(V) = 0$ and $\rho(V)^{-1}:=\infty$ in that case.
\end{theorem}
\begin{proof}
    Suppose first that $\rho(V)=0$. Then $V$ is nilpotent, and hence $(I-tV)^{-1}$ is a polynomial in $t$. As a consequence, $\Phi(A,t)$ is also polynomial in $t$, and therefore its radius of convergence is $\infty$ yielding the statement.

Now, assume that $\rho(V)>0$. Note that $\rho:=\rho(V)$ is the Perron eigenvalue of the non-negative matrix $V$. In particular, having fixed a vector norm $\| \cdot \|$ on $\R^m$, then there is an unique non-negative Perron right eigenvector $w \geq 0$ such that $V w = \rho w$ and $\| w \|=1$ \cite[Theorem 6.2.1]{Friedland}. Note that $\Phi(A,t)$ is the Schur complement of $tV-I$ in $P(t)$, where $P(t)$ is the pencil
\[  P(t)=\begin{bmatrix}
    tV - I & \sqrt{Z}R\\
    t L^T \sqrt{Z} & I
\end{bmatrix}.      \]
By \eqref{eq:genfunweighted} (see also \cite[Theorems 4.7 and 4.10]{rosenedd} for a connection to Rosenbrock's theorem \cite[Chapter 3, Theorem 4.1]{Rosen}), the finite poles of $\Phi(A,t)$ are always a subset of the finite eigenvalues of the pencil $tV-I$, which in turn are the inverses of the nonzero eigenvalues of $V$. Moreover, suppose that $\lambda^{-1}$ is a nonzero eigenvalue of $V$; then,  by \cite[Theorem 3.5]{local}, a sufficient condition for $\lambda$ to be a pole of $\Phi(A,t)$ is that $P(t)$ is locally minimal at $\lambda$  \cite[Definition 3.1]{local}. Equivalently, a sufficient condition is that the matrices
\[  \begin{bmatrix}
    \lambda V - I\\
    \lambda L^T \sqrt{Z}
\end{bmatrix}  \quad\text{and}\quad \begin{bmatrix}
    \lambda V-I & \sqrt{Z}R
\end{bmatrix}  \]
have, respectively, full column rank and full row rank.

Since the radius of convergence of a power series that converges to a rational function is the smallest pole of the rational function, we conclude that a sufficient condition for $\rho^{-1}$ to be the radius of convergence of $\Phi(A,t)$ is that $P(t)$ is locally minimal at $\rho^{-1}$. Suppose now for a contradiction that \[ \begin{bmatrix}
    \rho^{-1} V - I\\
    \rho^{-1} L^T \sqrt{Z}
\end{bmatrix} = \rho^{-1} \begin{bmatrix}
     V - \rho I\\
 L^T \sqrt{Z}
\end{bmatrix} \]
does not have full column rank. Then, there is a nonzero vector $x$ such that $(V-\rho I)x=0$ and $L^T \sqrt{Z} x=0$. The first equation is equivalent to $x$ being a Perron eigenvector of $V$, and hence without loss of generality we can assume that $x=w \geq 0$. The second equation says that the non-negative vector $y=\sqrt{Z} w$ belongs to $\ker L^T$. Recall that, by definition, $L_{ei}=1$ if and only if, in the underlying digraph $G$, edge $e$ starts from node $i$, and $L_{ei}=0$ otherwise. Fix now a node index $j$. As a consequence of the previous argument, $0=(L^T y)_j = \sum_e y_e$ where the sum is taken over all directed edges $e \in E(G)$ that start from node $j$. But $y_e \geq 0$, and hence $y_e=0$ for all such edges. Since $j$ is a generic node index, and since every edge start from some node in $G$, this arguments leads to $y=0$. Thus, $w=\sqrt{Z^{-1}}y=0$, which is a contradiction since $\|w\|=1$. A completely analogous argument, but starting with a Perron left eigenvector of  $V$, shows that 
\[\begin{bmatrix}
    \rho^{-1} V - I & \sqrt{Z} R
\end{bmatrix} = \rho^{-1} \begin{bmatrix}
    V - \rho I & \rho \sqrt{Z} R
\end{bmatrix}  \]
must have full row rank.

Consequently, $P(t)$ is locally minimal at $\rho^{-1}$ and thus $\rho^{-1}$ is the smallest pole of $\Phi(A,t)$. The statement immediately follows.
\end{proof}

 In Theorem \ref{th:radiusweight} below, we characterize the radius of convergence for the weighted case.
\begin{theorem} \label{th:radiusweight} 
Let $G$ be a weighted digraph. Assume that $G$ has $s$ strongly connected {or single node} components and denote by $H_1, \dots, H_s$ the undirectizations of each of them. Let $V$ be the matrix in \eqref{eq:Hashimoto_weight}, and let $\sigma$ be the smallest positive entry of $V$. Let $r$ be the radius of convergence of the generating function in \eqref{eq:genfunweighted}. Then the following statements hold:
\begin{itemize}
    \item[\rm(a)] If all the undirectizations $H_i$ are trees {or isolated nodes,} then $r=\infty$.
    \item[\rm(b)] If no one of the undirectizations $H_i$ has more than one cycle but at least one of them has precisely one cycle then $r\leq 1/\sigma.$
    \item[\rm(c)] If at least one undirectization $H_i$, for some $1 \leq i \leq s$, has more than one cycle then $r<1/\sigma.$
\end{itemize}
\end{theorem}

\begin{proof}
\begin{itemize}
   \item[\rm(a)] If all the undirectizations $H_i$ are trees or isolated nodes then $V$ is nilpotent, and hence $\rho(V)=0$. Then, by Theorem \ref{thm:weightedradius}, $r=\infty$.
    \item[\rm(b)] Consider the scaled matrix $\widehat V:=\frac{1}{\sigma}V$. Then, let $\widetilde G$ be the digraph obtained by setting all the weights of $G$ to be $1$, that is, the unweighted version of $G$. We now consider the Hashimoto matrix $\widetilde B$ of $\widetilde G$, and note that $\widehat V\geq \widetilde B$.  If no one of the undirectizations $H_i$ has more than one cycle but at least one of them has precisely one cycle then, by Theorem \ref{th:radius}$\rm(b)$ and Theorem \ref{thm:weightedradius}, $\rho(\widetilde B)=1$. Therefore, $\rho ( \widehat V)\geq \rho(\widetilde B) = 1$ and, thus, $\rho (V) \geq \sigma$. Then, by Theorem \ref{thm:weightedradius}, $r\leq 1/\sigma.$ 
    \item[\rm(c)] We use the same notation as in item (b). If at least one undirectization $H_i$, for some $1 \leq i \leq s$, has more than one cycle then by Theorem \ref{th:radius}$\rm(c)$ and Theorem \ref{thm:weightedradius}, $\rho(\widetilde B)>1$. Therefore, $\rho(\widehat V)\geq \rho(\widetilde B)>1$ and, thus, $\rho (V) > \sigma$. Then, by Theorem \ref{thm:weightedradius}, $r<1/\sigma.$
\end{itemize}
\end{proof}

\begin{remark}\rm
    Note that $\sigma$ in Corollary \ref{th:radiusweight} is equal to
\[  \min_{(e,f) \ \mathrm{is} \ \mathrm{a} \ NBTW} \sqrt{\omega(e)\omega(f)} .   \]
In other words, $\sigma$ is the minimum of the geometric mean of the weights of any pair of edges that form a NBTW of length $2$ in the graph.
\end{remark}

We finally prove a version of Ihara's theorem for weighted digraphs in Theorem \ref{th:det_weights} below. 

\begin{remark}\label{rem:notkempton} \rm
It is worth noting that another result appeared in the literature and was called ``weighted Ihara's Theorem'', namely, \cite[Theorem 4]{Kempton}. However, we emphasize that our setup for Theorem \ref{th:det_weights} is different and much more general. Indeed, (a) in \cite{Kempton} only weighted undirected graphs are considered, and (b) \cite{Kempton} focuses on the special case where each \emph{vertex} has a weight $\varphi(i)$ and the weight of the edge $e=(i,j)$ is $\omega(e)=\varphi(i)\varphi(j)$, whereas we allow the weight of each edge to be fully independent from each other.
\end{remark}

\begin{theorem}\label{th:det_weights}
Let $G$ be a weighted digraph and let $\Phi(A,t):=  I + t L^T \sqrt{Z} (I-tV)^{-1} \sqrt{Z} R$ be the generating function in \eqref{eq:genfunweighted}. Assume that $G$ has $b$ (undirected) reciprocal edges, and that each of them have weights $w(i)$ and $w'(i)$, respectively, for $i=1,\ldots,b$. Then

$$ \det \Phi(A,t)= \dfrac{\displaystyle \prod_{i=1}^{b} \left( 1- t^2 w(i)w'(i)  \right)}{\det (I-tV)}.$$

\end{theorem}

\begin{proof}
	Let $W_u:=RL^T$ be the line graph adjacency matrix of the unweighted version of $G$, i.e., the graph $\widetilde G$ having the same nodes and edges as $G$ but with all weights set equal to $1$. Moreover, let $\Delta$ be defined as in Section \ref{sec:BTDW}. Then, $B_u=W_u-\Delta$ is the Hashimoto matrix of $\widetilde G$ (see Section \ref{sec:BTDW} and \cite[Theorem 5.3]{nonBT}) while by \cite[Theorem 4.2]{ryan} $W=Z W_u Z$ is the adjacency matrix of the (weighted) line graph of $G$ (with the weight of an edge in the line graph being equal to the product of the weights of its two edges in $G$) and the weighted Hashimoto matrix of $G$ is $B=Z B_u Z$ (see Section \ref{sec:BTDW} and \cite{ryan}), since $Z$ is diagonal and positive definite and $B_u$ has entries all lying in $\{0,1\}$ we conclude that $V= \sqrt{Z} B_u \sqrt{Z}$.
 
 By the generalized matrix determinant lemma and \eqref{eq:genfunweighted}, we have that 
	\begin{equation}\label{eq:det_weight}
	\det (I-tV + t \sqrt{Z} W_u \sqrt{Z})=\det (I-tV) \det \Phi(A,t).
	\end{equation}
	On the other hand, taking into account that $V= \sqrt{Z} (W_u-\Delta) \sqrt{Z}$,
 \begin{equation}\label{eq:det_num}
	 \det (I-tV + t \sqrt{Z}W_u \sqrt{Z})= \det (I + t \sqrt{Z} \Delta \sqrt{Z} )= \det (I + t  \Delta Z ),
  \end{equation}
	since $ \sqrt{Z} $ is invertible. Then, by the proof of Lemma \ref{lem:det}, we have that 
	\begin{equation*}
	\Delta Z = \underbrace{ \begin{bmatrix}
		0 & w'(1)\\
		w(1) & 0
		\end{bmatrix} \oplus \dots \oplus \begin{bmatrix}
		0 & w'(b)\\
		w(b) & 0
		\end{bmatrix}}_{b \ \mathrm{copies}} \oplus \underbrace{0 \oplus \dots \oplus 0}_{a \ \mathrm{copies}}.
	\end{equation*}
	Hence, $\det (tI-\Delta Z)=t^a \prod_{i=1}^{b} \left( t^2-  w(i)w'(i)  \right)$. Thus, 
	\begin{align*}
	\det (I + t  \Delta Z ) & = (-t)^{a+2b} \det (-\frac{1}{t}I -   \Delta Z )=(-t)^{a+2b} \left( -\frac{1}{t} \right)^{a} \prod_{i=1}^{b} \left( \frac{1}{t^2}- w(i)w'(i)  \right)\\ & =  \prod_{i=1}^{b} \left( 1- t^2 w(i)w'(i)  \right).
	\end{align*}
	The statement follows by \eqref{eq:det_weight}, \eqref{eq:det_num} and the equation above.
\end{proof}
\begin{remark}\rm
In the classical case of unweighted digraphs, $V=B$ and $\omega(i)=1$ for all $i$ so Theorem \ref{th:det_weights} reduces to Corollary \ref{cor:Ihara}, and we have that $\det \Phi(A,t)^{-1}$ and $\det(I-tB)$ have the same finite roots with the same multiplicities except possibly $\pm 1$. In the weighted case, Theorem \ref{th:det_weights} shows that there are more possible exceptions, namely,  $\pm\dfrac{1}{\sqrt{\omega(i)\omega'(i)}}$.
\end{remark}

\section{Conclusions}\label{sec:conc}

We have studied the generating function associated to non-backtracking walk on a graph. We have extended the characterization of  the radius of convergence in \cite{GHN}, that related it to the number of cycles in the underlying graph in the unweighted and undirected case,  to possibly weighted and/or directed graphs. We showed that a similar characterization exists provided that one replace all the edges in the digraph by reciprocal edges. For weighted digraphs (or weighted undirected graphs), we showed for the first time that the radius of convergence is the inverse of the spectral radius of the weighed version of the Hashimoto matrix; this also extends a previously known result for the case of simple graphs. Finally, we have studied the case of backtrack-downweighted walks and we proved that a version of Ihara's theorem holds.


\newpage

\begin{thebibliography}{9}

\bibitem{indexsum} L.M. Anguas, F.M. Dopico, R. Hollister, and D.S. Mackey.
\newblock{\em Van Dooren's index sum theorem and rational matrices with prescribed structural data}.
\newblock SIAM J. Matrix Anal. Appl. 40, 720--738, 2019.

\bibitem{AGHN} F. Arrigo, P. Grindrod, D. J. Higham and V. Noferini.
\newblock {\em Non-backtracking walk centrality for directed networks}.
\newblock J. Complex Networks 6 (1), 54--78, 2018.

\bibitem{inprep} F. Arrigo, D. J. Higham and V. Noferini.
\newblock {\em  A theory for backtrack-downweighted walks}.
\newblock SIAM J. Matrix Anal. Appl. 42 (3), 1229--1247, 2021.

\bibitem{nonBT} F. Arrigo, D. J. Higham and V. Noferini,
\newblock	{\em Beyond non-backtracking: non-cycling network centrality measures.} \newblock Proc. R. Soc. A 476:20190653, 2020. 

\bibitem{ryan} F. Arrigo, D. J. Higham, V. Noferini and R. Wood. \newblock {\em Weighted enumeration of nonbacktracking walks on weighted graphs.} \newblock SIAM J. Matrix Anal. Appl., To appear, 2023.

\bibitem{Abel} H. Delange,
\newblock	{\em The converse of Abel's Theorem on power series.}
\newblock Ann. of Math. 50, 94--109, 1949.

\bibitem{local} F.~Dopico, S.~Marcaida, M.C.~Quintana, P.~Van Dooren,
\newblock {\em Local linearizations of rational matrices with application to rational approximations of nonlinear eigenvalue problems},
\newblock  Linear Algebra Appl. 604, 441--475, 2020.

\bibitem{rosenedd} F. Dopico, V. Noferini, I. Zaballa.
\newblock{\em Rosenbrock’s theorem on system matrices over
Elementary Divisor Domains}. 
\newblock In preparation.

\bibitem{flanders} H. Flanders.
{\em Elementary divisors of $AB$ and $BA$}.
Proc. Amer. Math. Soc. 2 (6), 871--874, 1951.

\bibitem{Friedland} S.\ Friedland,
\newblock {\it Matrices: Algebra, Analysis and Applications},
\newblock World Scientific, 2015.

\bibitem{GLR82} I.\ Gohberg, P.\ Lancaster, L.\ Rodman,
\newblock {\it Matrix Polynomials},
\newblock Academic Press, 1982.

\bibitem{GHN} P. Grindrod, D. J. Higham and V. Noferini.
\newblock {\em The deformed graph Laplacian and its applications to network centrality analysis}.
\newblock SIAM J. Matrix Anal. Appl. 39 (1), 310--341, 2018.

\bibitem{GT} S. G\"{u}ttel, F. Tisseur, {\em The nonlinear eigenvalue problem},  Acta Numer. {\bf 26}: 1–94 (2017).

\bibitem{Hashimoto} K. Hashimoto, 
\newblock {\em On zeta and L-functions of finite graphs}. 
\newblock Int. J. Math. 1, 381--396, 1990.


\bibitem{HST06}
M. D. Horton,  H. M. Stark,  and A. A. Terras. {\em What  are  zeta  functions  of  graphs  and what  are  they  good  for?},  in: Quantum  Graphs  and  Their  Applications,  Contemp.  Math.415, G. Berkolaiko, R. Carlson, S. A. Fulling, and P. Kuchment, eds., 2006, pp. 173–190.

\bibitem{ihara} Y. Ihara. {\em On discrete subgroups of the two by two projective linear group over p-adic fields}. J. Math.
Soc. Japan 18, 219--235,1966.

\bibitem{Katz} L. Katz, {\em A new status index derived from sociometric data analysis}, Psychometrika 18, 39--43, 1953.

\bibitem{Kempton} M. Kempton, {\em Non-Backtracking Random Walks and a Weighted Ihara’s Theorem}, Open Journal of Discrete Mathematics, 6, 207--226, 2016.


\bibitem{Morbidi} F. Morbidi. \newblock {\em The deformed consensus protocol}.
\newblock Automatica 49, 3049--3055, 2013.

\bibitem{mulas} R. Mulas, D. Zhang, G. Zucal. \newblock {\em There is no going back: Properties of the non-backtracking Laplacian.} \newblock Preprint, \url{https://arxiv.org/pdf/2303.00373.pdf}

\bibitem{nnt} Y.\ Nakatsukasa, V.\ Noferini, A.\ Townsend,
\newblock {\em Vector spaces of linearizations for matrix polynomials: a bivariate polynomial approach}, 
\newblock  SIAM J. Matrix Anal. Appl.  38(1), 1--29, 2017.

\bibitem{pc20} R. Pastor-Satorras and C. Castellano. \newblock {\em The localization of non-backtracking centrality in networks and its physical consequences.} Sci. Rep 10., 21639, 2020.

\bibitem{Rosen} H. H. Rosenbrock,
\textit{State-space and Multivariable Theory,} Thomas Nelson and Sons, London, 1970.

\bibitem{tarfulea} A. Tarfulea and R. Perlis.
\newblock {\em An Ihara formula for partially directed graphs}.
\newblock Linear Algebra Appl. 431, 73--85, 2009.

\bibitem{tcter} L. Torres, K. S. Chan, H. Tong and T. Eliassi-Rad. \newblock {\em Nonbacktracking eigenvalues under node removal: X-Centrality and targeted immunization.} \newblock SIAM J. Math. Data Sci. 3(2), 656--675, 2021.






\end{thebibliography}
\end{document}